\documentclass[hidelinks,onefignum,onetabnum]{siamart220329}

\usepackage{enumerate}
\usepackage{cite}

\usepackage{amssymb,amsmath,bm,euscript}
\usepackage{graphicx,color,caption}
\usepackage{ifpdf}
\usepackage{epstopdf}
\usepackage{bm}

\usepackage{braket}
\usepackage{setspace}
\usepackage{algorithm, algorithmic}
\usepackage{multirow}
\usepackage{multicol}
\usepackage{subfigure}
\usepackage{booktabs}
\usepackage[figuresright]{rotating}

\usepackage[OT1]{fontenc}

\def\A{{\mathcal{A}}}
\def\B{\mathcal{B}}
\def\C{\mathcal{C}}
\def\CC{\mathbb{C}}
\def\KK{\mathbb{K}}
\def\RR{\mathbb{R}}

\def\X{{\mathcal{X}}}
\def\Y{{\mathcal{Y}}}

\def\aa{{\bf a}}
\def\bb{{\bf b}}
\def\cc{{\bf c}}

\def\ii{{\bf i}}

\def\ee{{\bf e}}

\def\0{{\bf 0}}
\def\1{{\bf 1}}

\usepackage{epstopdf}

\usepackage{lipsum}
\usepackage{amsfonts}
\usepackage{graphicx}
\usepackage{epstopdf}
\usepackage{algorithmic}
\ifpdf
  \DeclareGraphicsExtensions{.eps,.pdf,.png,.jpg}
\else
  \DeclareGraphicsExtensions{.eps}
\fi

\newsiamremark{remark}{Remark}
\newsiamremark{hypothesis}{Hypothesis}
\crefname{hypothesis}{Hypothesis}{Hypotheses}
\newsiamthm{claim}{Claim}

\headers{Variable T-Product and Tensor Completion}{Liqun Qi, Rui Yan, Ziyan Luo, Hong Yan, and Gaohang Yu}
\title{Variable T-Product and Zero-Padding Tensor Completion with Applications\thanks{Submitted to the editors DATE.
\funding{This work was supported in part by National Natural Science Foundation of China (No.12071104), Natural Science Foundation of Zhejiang Province (No.LD19A010002), Beijing Natural Science Foundation (No.Z190002), Hong Kong Research Grants Council (Project 11204821), Hong Kong Innovation and Technology Commission(InnoHK Project CIMDA) and City University of Hong Kong (Project 9610034).}}}

		\author{Liqun Qi\thanks{Department of Mathematics, Hangzhou Dianzi University, Hangzhou, 310018, China 
		(\email{liqun.qi@polyu.edu.hk}, zhxsryan@hotmail.com, maghyu@hdu.edu.cn).}
	\and Rui Yan\footnotemark[2]
	\and Ziyan Luo\thanks{School of Mathematics and Statistics,
		Beijing Jiaotong University, Beijing 100044, China
		(\email{starkeynature@hotmail.com}).}
	\and Hong Yan\thanks{Department of Electrical Engineering, and Center for Intelligent Multidimensional Data Analysis, City University of Hong Kong, Kowloon, Hong Kong, China
		(\email{h.yan@cityu.edu.hk}).}
	\and Gaohang Yu\footnotemark[2]}

\usepackage{amsopn}

\ifpdf
\hypersetup{
  pdftitle={An Example Article},
  pdfauthor={D. Doe, P. T. Frank, and J. E. Smith}
}
\fi

\begin{document}

\maketitle

\begin{abstract}
The T-product method based upon Discrete Fourier Transformation (DFT) has found wide applications in engineering, in particular, in image processing.
In this paper, we propose variable T-product, and apply the Zero-Padding Discrete Fourier Transformation (ZDFT), to convert third order tensor problems to the variable Fourier domain.  An additional positive integer parameter is introduced, which is greater than the tubal dimension of the third order tensor.   When the additional parameter is equal to the tubal dimension, the ZDFT reduces to DFT. Then, we propose a new tensor completion method based on ZDFT and TV regularization, called VTCTF-TV. Extensive numerical experiment results on visual data demonstrate that the superior performance of the proposed method.    In addition,  the zero-padding is near the size of the original tubal dimension, the resulting tensor completion can be performed much better. 
\end{abstract}



\section{Introduction}
The Discrete Fourier Transformation (DFT) based the tensor SVD (T-SVD) factorizations were proposed by Kilmer, Martin and others \cite{KM11, KBHH13}, and have found extensive applications for tensor completion in machine learning, imaging processing, quantum information and other engineering fields \cite{kolda2009tensor, CXZ20, SHKM14, WGLZ21, XCGZ21, xiao2020prior, YHHH16, ZSKA18, ZA17, ZEAHK14, ZLLZ18,yu2022t,shi2021robust,tarzanagh2018fast,zhao2022robust}.   In 2015, Knerfeld, Kilmer and Aeron \cite{KKA15} further extended T-SVD to the so-called transformed T-SVD by replacing the DFT with an arbitrary invertible linear transformation.   Since then, transformed T-SVD has attracted attention of many researcheers \cite{HZW21, JNZH20, qiu2021robust, qiu2021nonlocal, SNZ20, ZC21,li2022nonlinear,deng2023t,zheng2020tensor,zheng2019mixed}.  

A typical application of the T-product method is for color image processing, in which the data tensor is an $m \times n \times p$ real third order tensor $\A = (a_{ijk})$.   Here, $ m \times n $ is the image size, and $ p=3 $ for three color channels, red, green, and blue (RGB).  One way to handle this problem is to use a quaternion matrix to represent the third order tensor $\A$ \cite{chen2022color, LS03, MKL20, qi2022quaternion}.   Another way to deal with this problem is to regard $\A$ as a tubal matrix $\A = (\aa_{ij})$, where the entries of the tubal matrix $\A$ are $p$-dimensional vectors.    For each pixel, a $p$-dimensional vector of data is acquired.  Then, via the DFT, the scalar product of entries in the classical matrix product is replaced by a convolution-like product of the vector entries of the tubal matrix.   In this regard, the matrix SVD can be extended to T-SVD.  This is the essence of the T-product method.

Since the T-SVD method suffers from expensive computational cost, Zhou, Lu, Lin and Zhang \cite{ZLLZ18} proposed a tensor factorization model for the low rank tensor completion problem (TCTF).    Instead of analyzing the data in the original image domain, they solve an optimization problem in the Fourier domain.   If the original data tensor is in $\mathbb{R}^{m \times n \times p}$, then the Fourier domain is in $\mathbb{C}^{m  \times n \times p}$. 

However, in signal processing and image processing, direct application of the DFT may not be enough.  Zero-Padding \cite{MWGdD02, DT15} is needed to obtain the accurate result.   Otherwise, artifacts may be produced. With the above motivation, in this paper, we propose variable T-product, and apply the Zero-Padding Discrete Fourier Transform (ZDFT) to the data tensor.   Then we solve an optimization problem in the variable Fourier domain in $\mathbb{C}^{m \times  n \times v}$.   Here, $v$ is an additional positive integer parameter, with $v > p$. If $v=p$, then the  variable Fourier domain and the ZDFT reduce to the Fourier domain and the DFT. The ZDFT matrix consists of the first $p$ columns of the $v \times v$ DFT matrix, and it maps an original $m \times n \times p$ tensor  to the $m \times n \times v$ variable Fourier domain. In theory, shorter signals $(v< 2p -1)$ would produce artifacts, and longer zero-padding $(v> 2p-1)$ should produce the same result, but the result can be slightly different due to quantization of the frequency variable in the Fourier domain. Therefore, when $v = 2p-1$ or near $2p$, the Zero-Padding takes effect, and the resulting tensor completion can be performed much better.  Numerical experiments show that the new method performs better when $v$ is $2p-1$.

Throughout this paper, we assume that $m, n, p$ and $v$ are positive integers, and $v > p$.
\subsection{Contributions}
Our contributions are summarized as below:
\begin{itemize}
	\item[(1)] We propose variable T-product and apply the Zero-Padding Discrete Fourier Transform (ZDFT), to convert third order tensor problems to the variable Fourier domain $\mathbb{C}^{m \times n \times v}$.
	\item[(2)] We present a novel Zero-Padding tensor completion method, termed as VTCT\ F-TV, for low-rank tensor completion problem, based on the ZDFT. Numerical experiments verify the superior performance of the proposed method in image and video inpainting when $v$ is or near $2p-1$.
	\item[(3)] We design an efficient proximal alternating minimization
	(PAM) algorithm to solve the proposed tensor completion problem and establish its convergence results.
\end{itemize}

\section{Preliminaries}

\subsection{Tensor Notations and Definitions}
An $M$th-order tensor $\mathcal{X}\!\in\!\mathbb{R}^{I_1\times\dots\times I_M}$ is an $M$-way array consisting of entries $x_{i_1 \cdots i_m}$ with each $i_j$ varying among ${1,\cdots,I_j}$ for all $j=1,\cdots, M$. Vectors and matrices are typical low order tensors with $M=1$ and $M=2$, respectively. 
Define the inner product, the Frobenius norm and the $\infty$ norm for tensors with symbols $\langle\cdot,\cdot\rangle$, $\|\cdot\|_{F}$, $\|\cdot\|_{\infty}$ respectively, formulated by
$$\langle\mathcal{X}, \mathcal{Y}\rangle:=\sum_{i_1 = 1}^{I_1}\dots\sum_{i_M = 1}^{I_M} x_{i_1\dots i_M} y_{i_1\dots i_M},$$
$$\Vert\mathcal{X}\Vert_F:=\sqrt{\langle\mathcal{X}, \mathcal{X}\rangle},~~ \|\mathcal{X}\|_{\infty}: = \max\limits_{i_1, \cdots, i_M} | x_{i_1\dots i_M}|,$$
for all $\mathcal{X}, \mathcal{Y}\in\mathbb{R}^{I_1\times\dots\times I_M}$.

The Hadamard product of two vectors or matrices or tensors with the same dimensions, is the component-wise product of these two vectors, or matrices, or tensors.

In this paper, we deal with third order real tensors. A third order tensor $\A = (a_{ijk}) \in \mathbb{R}^{m \times n \times p}$ has three modes:  row, column and tubal.  Let $A^{(k)} = \A(:, :, k)$ be the $k$th frontal slice of $\A$ for $k = 1, \cdots, p$. Given a positive integer $t\leq p$, let $\A(:, :, [1:t])\in \RR^{m\times n\times t}$ be the subtensor of $\A$, generated by the first $t$ frontal slices.   As stated in the introduction, we may regard the third order tensor $\A = (a_{ijk})$ as a tubal matrix $\A = (\aa_{ij})$, where $\aa_{ij} = \A(i,j,:)$ is a $p$-dimensional vector, also called a tubal scalar \cite{KBHH13}.

We use calligraphic letters $\A, \B, \C, \cdots$ to denote tensors and tubal matrices, capital letters $A, B, C, \cdots$ to denote matrices, small bold letters $\aa, \bb, \cc$ to denote vectors and tubal scalars, small letters $a, b, c. \cdots$ to denote scalars.

\subsection{The DFT process}\label{SubsectionIIB}

Let $\C \in \mathbb{R}^{m \times n \times p}$ be the data tensor.   Let $F = F_p$ be the $p \times p$  Discrete Fourier Transform (DFT) matrix with the form
$$F =  \begin{bmatrix}
1 & 1 & \cdots & 1 & 1 \\
1 & \omega_p & \cdots  & \omega_p^{p-2} & \omega_p^{p-1} \\
\vdots & \vdots & \ddots& \vdots & \vdots \\
1 & \omega_p^{p-2} &  \cdots & \omega_p^{(p-2)(p-2)} & \omega_p^{(p-2)(p-1)}\\
1 & \omega_p^{p-1} &  \cdots & \omega_p^{(p-1)(p-2)} & \omega_p^{(p-1)(p-1)}
\end{bmatrix},$$
where
$$\omega_p = e^{-{2\pi \ii \over p}},$$
and $\ii = \sqrt{-1}$ is the imaginary unit.   Then $F$ is invertible and $F^{-1} = {1 \over p}F^H$.   Apply $F$ to $\C$ along the third dimension, we have
$\bar \C(F) \in \mathbb{C}^{m \times n \times p}$ such that
\begin{equation} \label{ee1}
\bar \C(F)(i,j,:) = F\C(i,j, :),
\end{equation}
for any $i$ and $j$.  Then
\begin{equation} \label{ee2}
\C(i,j,:) = {1 \over p}F^H \bar \C(F)(i,j, :),
\end{equation}
for any $i$ and $j$, and
\begin{equation}
\|\C\|_F^2 = {1 \over p}\|\bar \C(F) \|_F^2.
\end{equation}
\par Define $\bar C(F) \in \mathbb{C}^{mp \times np}$ as a block diagonal matrix
\begin{equation}
\bar C(F) = {\rm diag}\left(\bar C(F)^{(1)}, \bar C(F)^{(2)}, \cdots, \bar C(F)^{(p)}\right).
\end{equation}

\begin{definition}
	For $\C \in \mathbb{R}^{m \times n \times p}$, its tensor tubal rank is defined as
	$${\rm Rank}_p(\C) = \max \{ {\rm Rank}(\bar C(F)^{(1)}), \cdots, {\rm Rank}(\bar C(F)^{(p)})\}.$$
\end{definition}

\subsection{The TCTF Method}

Since T-SVD suffers from expensive computational cost, Zhou, Lu, Lin and Zhang \cite{ZLLZ18} proposed a tensor factorization model for the low rank tensor completion problem (TCTF). They work in the Fourier domain by solving the following optimization problem:
\begin{equation}\label{opt:TCTF}
\begin{aligned}
\mathop {\min }\limits_{\hat {\cal X}, \hat {\cal Y},\C} \quad & \frac{1}{2p} \sum_{k=1}^p \left\| {\hat {X}^{(k)} \hat {Y}^{(k)} - \bar C(F)^{(k)}} \right\|_F^2\\
\mbox{s.t.} \quad  & {{\cal P}_\Omega }({\cal C} - {\cal G}) = 0\\
&\hat \X \in \mathbb{C}^{m\times r\times p}, \hat \Y \in \mathbb{C}^{r\times n \times p},\mathcal{C}\in \mathbb{R}^{m\times n\times p},
\end{aligned}
\end{equation}
where $\C$ and $\bar \C(F)$ are related by (\ref{ee1}) and (\ref{ee2}), $\mathcal{G}=(g_{i_1i_2i_3})$ is the observed incomplete tensor data in $\mathbb{R}^{m\times n\times p}$, $\Omega$ is the index set corresponding to the observed  entries of $\mathcal{G}$, and $P_{\Omega}$ is the linear operator to extract known elements in the subset $\Omega$ and to fill the elements that are not in $\Omega$ with zero values. Here, $r$ is the beforehand estimated tubal rank and is usually much smaller than ${\min}\{m,n\}$.

\section{Main Results}

In Subsection \ref{SubsectionIIIA}, we provide the definition of variable T-product of two third order tensors.  In Subsection \ref{SubsectionIIIB}, we describe the ZDFT.   We define the variable tensor tubal rank in Subsection \ref{SubsectionIIIC}.  A Zero-Padding tensor completion method is proposed in Subsection \ref{SubsectionIIID}.

\subsection{Variable T-Product}\label{SubsectionIIIA}

We first define the variable product between two vectors in $\RR^p$.  For $\aa \in \RR^p$ or $\CC^p$, denote its $i$th component as $\aa(i)$.

\begin{definition}[Variable product of two $p$-dimensional vectors]
	For any $\aa, \bb \in \RR^p$ or $\CC^p$, we have
	\begin{eqnarray} \label{e2.1}
	&&(\aa \odot_v \bb)(k)\\
	&=&\sum \left\{ \aa(i)\bb(j) : \begin{array}{l}
	i+j-k-1= 0\ {\rm mod}(v), \\
	i, j =1, \cdots, p
	\end{array}
	\right\}, \nonumber
	\end{eqnarray}
	for $k = 1, \cdots, p$.   We call $\aa \odot_v \bb$ the variable product of $\aa$ and $\bb$.
\end{definition}

One can see that when $v= 2p-1$, $\aa \odot_v \bb$ is exactly the $p$-truncated convolution of $\aa$ and $\bb$, since it consists of the first $p$ components of the convolution. It is well-known that the convolution of two $p$-dimensional vectors is a $(2p-1)$-dimensional vector.  In order that the variable T-product of two third order tensors can be defined, we need to keep the variable product of two $p$-dimensional vectors in dimension $p$.  Thus, we cannot use convolution of two vectors to define the variable product of two vectors directly here.  But convolution is the best way to combine two $p$-dimensional signals, and the variable product is the closest to the convolution when $v=2p-1$. This explains why the tensor completion method proposed later performs much better when $v=2p-1$.

Let $\ee \in \RR^p$ be defined by $\ee(1) = 1$, $\ee(k) = 0$, for $k = 2, \cdots, p$.

\begin{proposition}
	For any $\aa, \bb \in \RR^p$, we have
	\begin{equation} \label{e2.2}
	\aa \odot_v \bb = \bb \odot_v \aa,
	\end{equation}
	and
	\begin{equation} \label{e2.3}
	\aa \odot_v \ee = \ee \odot_v \aa = \aa.
	\end{equation}
\end{proposition}
\begin{proof} 
By (\ref{e2.1}), we have (\ref{e2.2}).  For $k = 1, \cdots, p$, we have
\begin{eqnarray*}
	&&(\aa \odot_v \ee)(k) \nonumber\\
	& = & \sum \left\{ \aa(i)\ee(j) : \begin{array}{l}
		i+j-k-1= 0\ {\rm mod}(v), \\
		i, j =1, \cdots, p
	\end{array} \right\}\\
	& = & \sum \left\{ \aa(i) : i= k, i = 1, \cdots, p \right\}\\
	& = & \aa(k).
\end{eqnarray*}
Thus, $\aa \odot_v \ee = \aa$.   Similarly, we have $\ee \odot_v \aa = \aa$.   This proves (\ref{e2.3}).
\end{proof}

Then $\KK_p = (\RR^p, +, \odot_v)$ is a commutative ring with unity $\ee$, with $+$ as the usual addition of vectors.  As discussed in \cite{Br10},  $\KK_p$ is also a module.   As in \cite{KBHH13}, we call an element $\aa \in \KK_p$ a tubal scalar, and call $\KK_p$ the $p$-dimensional tubal scalar module.
The zero tubal scalar $\0 \in \KK_p$ has components $\0(k) = 0$ for $k = 1, \cdots, p$.

For $\aa \in \KK_p$, define its modulus as
$$|\aa| = \sqrt{\sum_{k=1}^p \aa(k)^2}.$$

\begin{definition}
	Let $\aa \in \KK_p$.   We say that $\bb \in \KK_p$ is the transpose of $\aa$, and denote $\aa^\top = \bb$, if $\bb(1) = \aa(1)$, and
	$\bb(k)=\aa(p+2-k)$ for $k = 2, \cdots, p$.
\end{definition}

For any $\aa \in \KK_p$,
if $\aa^\top = \aa$, then we say that $\aa$ is symmetric. The following proposition can be proved by definition.

\begin{proposition} \label{p2.4}
	Let $\aa, \bb \in \KK_p$.  Then $(\aa^\top)^\top = \aa$ and $(\aa \odot_v \bb)^\top = \bb^\top \odot_v \aa^\top$.
\end{proposition}

As stated early, a third order tensor $\A = (a_{ijk}) \in \RR^{m \times n \times p}$ can be regarded as a tubal matrix $\A = (\aa_{ij})$, where $\aa_{ij}(k) = a_{ijk}$ for $k = 1, \cdots, p$.

\begin{definition}[Variable T-product of third order tensors in the real domain]
	Suppose that $\A = (\aa_{il}) \in \RR^{m \times q \times p}$ and $\B = (\bb_{lj}) \in \RR^{q \times n \times p}$.   The variable T-product of $\A$ and $\B$, denoted as $\C = \A *_v \B$,  is defined by $\C = (\cc_{ij})$ with
	\begin{equation} \label{e2.4}
	\cc_{ij} = \sum_{l=1}^q \aa_{il} \odot_v \bb_{lj}.
	\end{equation}
\end{definition} Apparently, if $v=p$, then ``$*_v$" reduces to the $T$-product ``$*$".

\subsection{Zero-Padding Discrete Fourier Transform}\label{SubsectionIIIB}

Let $T = T_{vp}$ be the $v \times p$  Zero-Padding Discrete Fourier Transform (ZDFT) matrix with the form
$$T =  \begin{bmatrix}
1 & 1 & \cdots & 1 & 1 \\
1 & \omega_v & \cdots  & \omega_v^{p-2} & \omega_v^{p-1} \\
\vdots & \vdots & \ddots& \vdots & \vdots \\
1 & \omega_v^{v-2} &  \cdots & \omega_v^{(v-2)(p-2)} & \omega_v^{(v-2)(p-1)}\\
1 & \omega_v^{v-1} &  \cdots & \omega_v^{(v-1)(p-2)} & \omega_v^{(v-1)(p-1)}
\end{bmatrix},$$
where
$$\omega_v = e^{-{2\pi \ii \over v}},$$
and $\ii = \sqrt{-1}$ is the imaginary unit.
Then $T$ consists of the first $p$ columns of the $v \times v$ DFT matrix $F_v$. Thus $T$ is of full column rank $p$, and $T^HT = vI_p$, where $I_p$ is the $p \times p$ identity matrix.  For any $\aa \in \CC^p$, zero-padding $\aa$ to obtain $\tilde{\aa}:=(\aa; {\bf 0})\in \CC^{v}$, we have \begin{equation}\label{zp}
T \aa = F_v \tilde{\aa}.\end{equation} For this reason, $T$ is called the ZDFT matrix.

We may see that the variable T-product can be interpreted by the ZDFT.

For $\aa \in \CC^p$ and $\cc \in \CC^v$,
define a mapping $\phi: \CC^p \to \CC^v$ by $\phi(\aa) := T\aa$, and its conjugate transpose mapping  $\phi^H : \CC^v \to \CC^p$ by $\phi^H(\cc) := T^H\cc$. Then for any $\aa \in \CC^p$ and $\cc \in \CC^v$, we have
\begin{equation} \label{e3.1}
\phi(\aa)(l) = \sum_{i=1}^p \omega^{(i-1)(l-1)} \aa(i),
\end{equation}
for $l=1, \cdots, v$, and
\begin{equation} \label{e3.2}
\phi^H(\cc)(k) = \sum_{l=1}^v \bar \omega^{(l-1)(k-1)} \cc(l),
\end{equation}
for $k = 1, \cdots, p$.

\begin{proposition} \label{p3.1}
	For any $\aa, \bb \in \RR^p$, we have
	\begin{equation} \label{e3.3}
	\aa \odot_v \bb = {1 \over v}\phi^H[\phi(\aa) \circ \phi(\bb)],
	\end{equation}
	where $\circ$ represents the Hadamard product.
\end{proposition}

\begin{proposition} \label{truncated-T}
	Suppose that $\A = (\aa_{il}) \in \RR^{m \times q \times p}$ and $\B = (\bb_{lj}) \in \RR^{q \times n \times p}$. Let $\A_0\in \RR^{m \times q \times v}$ and $\B_0\in \RR^{q \times n \times v}$ be the zero-padding counterparts of $\A$ and $\B$, with $\A_0(:,:,k)$ and $\B_0(:,:,k)$ all zero frontal slices for any $k= p+1, \cdots, v$. Then
	\begin{equation} \label{T-vT}
	\A *_v \B = \left(\A_0 * \B_0\right)(:,:,[1:p]),
	\end{equation}
	where $*$ is the T-product.
\end{proposition}

\subsection{Variable Tensor Tubal Rank and H-product}\label{SubsectionIIIC}
We now use $T$ to replace $F$ in Subsection \ref{SubsectionIIB}.  Again, let $\C \in \mathbb{R}^{m \times n \times p}$ be the data tensor.  Apply $T$ to $\C$ along the third dimension, we have
$\bar \C(T) \in \mathbb{C}^{m \times n \times v}$ such that
\begin{equation} \label{ee6}
\bar \C(T)(i,j,:) = T\C(i,j, :),
\end{equation}
for any $i$ and $j$.  Then
\begin{equation} \label{ee7}
\C(i,j,:) = {1 \over v}T^H \bar \C(T)(i,j, :),
\end{equation}
for any $i$ and $j$, and
\begin{equation}
\|\C\|_F^2 = {1 \over v}\|\bar \C(T) \|_F^2.
\end{equation}
Define $\bar C(T) \in \mathbb{C}^{mv \times nv}$ as a block diagonal matrix
\begin{equation}
\bar C(T) = {\rm diag}\left(\bar C(T)^{(1)}, \bar C(T)^{(2)}, \cdots, \bar C(T)^{(v)}\right).
\end{equation}

Clearly, the equations (\ref{ee6}) and (\ref{ee7}) define two mappings between $\C$ and $\bar C(T)$, respectively. For simplicity, let's use $\bar C(T)=\mbox{fft}_{p2v}(\C,[],3)$ to denote the relation (\ref{ee6}), and $\C=\mbox{ifft}_{v2p}(\bar C(T),[],3)$ to denote the relation (\ref{ee7}).

\begin{definition}[Variable tensor tubal rank]\label{v-tubal-rank}
	For $\C \in \mathbb{R}^{m \times n \times p}$, its variable tensor tubal rank is defined as
	$${\rm Rank}_v(\C) = \max \{ {\rm Rank}(\bar C(T)^{(1)}), \cdots, {\rm Rank}(\bar C(T)^{(v)})\}.$$
\end{definition} Obviously, if $v =p$, then ${\rm Rank}_v(\C)$ is exactly the tubal rank ${\rm Rank}_t(\C)$.

We may define H-product for tensors in $\mathbb{C}^{m \times n \times v}$.   Again, we regard tensors in $\mathbb{C}^{m \times n \times v}$ as tubal matrices.

\begin{definition} [H-product for third order tensors in the variable Fourier domain]
	Suppose that $\A = (\aa_{il}) \in \CC^{m \times q \times v}$ and $\B = (\bb_{lj}) \in \CC^{q \times n \times v}$.   The Hadamard tensor product, or simply called the H-product of $\A$ and $\B$, denoted as $\C = \A *_H \B$,  is defined by $\C = (\cc_{ij})$ with
	\begin{equation} \label{e2.24}
	\cc_{ij} = \sum_{l=1}^q \aa_{il} \circ \bb_{lj}.
	\end{equation}
\end{definition}

\begin{theorem} \label{t10}
	Suppose that $\C \in \RR^{m \times n \times p}$.   Then the variable tensor tubal rank of $\C$ is not greater than $r$ if and only if there are tensors $\bar{\X} \in \CC^{m \times r \times v}$ and $\bar{\Y} \in \CC^{r \times n \times v}$ such that $\bar \C(T) = \bar{\X} *_H \bar{\Y}$.   Furthermore, if there are $\A \in \RR^{m \times r \times p}$ and $\B \in \RR^{r \times n \times p}$ such that $\C = \A *_v \B$, then the  variable tensor tubal rank of $\C$ is not greater than $r$.
\end{theorem}

\subsection{A Zero-Padding Tensor Completion Method}\label{SubsectionIIID}

Based on the Zero-Padding Discrete Fourier Transform, we propose the following Zero-Padding Tensor Completion Method. In order to accurately recover the missing data, spatio-temporal prior knowledge is introduced, and TV regularization is applied to tensor completion problem by the application of total variation smoothing prior method in spatio-temporal video \cite{lin2020tensor}. In this paper, the low-rank completion problem of third-order tensor based on variable tubal rank and TV regularization is studied. The model is as follows:
\begin{equation}\label{eq:1}
\begin{aligned}
&\min_{\mathcal{X},\mathcal{Y},\mathcal{C}  } \frac{1}{2} \left \| \mathcal{X}\ast _{v}\mathcal{Y}-\mathcal{C}    \right \|_{F}^{2}+\alpha _{1}\left \| \mathcal{D}_{1}\ast _{v}\mathcal{C}     \right \|_{l_{1} }+\alpha _{2}\left \| \mathcal{C}\ast _{v} \mathcal{D} _{2}   \right \| _{l_{1} }\\
&s.t. \mathcal{C}_{\Omega }=\mathcal{G}_{\Omega }\\ 
&\mathcal{X}\in \mathbb{R}^{m\times q\times p}, \mathcal{Y}\in \mathbb{R}^{q\times n\times p},\mathcal{C}\in \mathbb{R}^{m\times n\times p},
\end{aligned}
\end{equation}
where $ \mathcal{G} $ represents the observation tensor, $ \Omega $ represents the index set of observation elements, and $ \mathcal{C}_{\Omega }=\mathcal{G}_{\Omega } $ represents that the elements of $ \mathcal{C} $ and $ \mathcal{G} $ in $ \Omega $ are consistent, and	
$$
\mathcal{D}_{1}\left ( :,:,i \right )=\begin{cases} L_{m},i=1 \\ O_{m},i> 1  \end{cases},
$$
$$
\mathcal{D}_{2}\left ( :,:,i \right )=\begin{cases} L_{n}^{\ast } ,i=1 \\ O_{n},i> 1  \end{cases},
$$
$$
L_{m}=\begin{pmatrix}
0 & 0 & 0 & \cdots  & 0 & 0\\
-1 & 1 & 0 & \cdots  & 0 & 0\\
0 & -1 & 1 & \cdots  & 0 & 0\\
\vdots  & \vdots  & \vdots  & \ddots  & \vdots  & \vdots \\
0 & 0 & 0 & \cdots  & 1 & 0\\
0 & 0 & 0 & \cdots  & -1 & 1
\end{pmatrix} \in \mathbb{R}^{m\times m}   , 	
$$
$ \alpha _{1} $ and $ \alpha _{2} $ are regularization parameters and $ \alpha _{1}\left \| \mathcal{D}_{1}\ast _{v}\mathcal{C}     \right \|_{l_{1} }+\alpha _{2}\left \| \mathcal{C}\ast _{v} \mathcal{D} _{2}   \right \| _{l_{1} } $ is TV regularization term. The above optimization problems can be transformed into the following unconstrained optimization problems:	
\begin{equation}\label{eq:2}
\min f\left ( \mathcal{X},\mathcal{Y},\mathcal{C}    \right ) =\frac{1}{2} \left \| \mathcal{X}\ast _{v}\mathcal{Y}-\mathcal{C}    \right \|_{F}^{2}+\alpha _{1}\left \| \mathcal{D}_{1}\ast _{v}\mathcal{C}     \right \|_{l_{1} }+\alpha _{2}\left \| \mathcal{C}\ast _{v} \mathcal{D} _{2}   \right \| _{l_{1} }+\Phi \left ( \mathcal{C}  \right ) ,
\end{equation}
where
$$
\Phi \left ( \mathcal{C}  \right ) =\begin{cases} 0,\mathcal{C}_{\Omega }= \mathcal{G}_{\Omega }, \\ \infty ,otherwise.  \end{cases} 
$$
\par Optimization problem (\ref{eq:2}) is not a joint convex function about $ \left ( \mathcal{X} ,\mathcal{Y} ,\mathcal{C}  \right )  $, but it is convex for each variable $ \mathcal{X} $, $ \mathcal{Y} $ and $ \mathcal{C} $. Because of the high efficiency of Alternate Minimization (AM) algorithm, it is often used to solve multivariate optimization problems. In order to improve the theoretical convergence and numerical stability of the algorithm, this paper adds a neighboring term to the sub-problem of AM algorithm, that is, the Proximal Alternate Minimization (PAM) algorithm.	
Let $ f\left ( \mathcal{X},\mathcal{Y},\mathcal{C}    \right ) $ be the objective function of problem (\ref{eq:2}) and given the initial point $ \left ( \mathcal{X}^{k},\mathcal{Y}^{k},\mathcal{C}^{k}   \right )  $, use PAM algorithm to solve the framework, and update each variable alternately as follows:
\begin{equation}\label{eq:33}
\begin{cases} \mathcal{X}^{k+1}=\underset{\mathcal{X}}{argmin }f\left ( \mathcal{X},\mathcal{Y}^{k},\mathcal{C}^{k}  \right )+\frac{\rho _{1} }{2}\left \| \mathcal{X}-\mathcal{X}^{k}    \right \| _{F}^{2}      \\ \mathcal{Y}^{k+1}=\underset{\mathcal{Y}}{argmin }f\left ( \mathcal{X}^{k+1},\mathcal{Y},\mathcal{C}^{k}  \right )+\frac{\rho _{2} }{2}\left \| \mathcal{Y}-\mathcal{Y}^{k}    \right \| _{F}^{2}\\ \mathcal{C}^{k+1}=\underset{\mathcal{C}}{argmin }f\left ( \mathcal{X}^{k+1},\mathcal{Y}^{k+1},\mathcal{C}  \right )+\frac{\rho _{3} }{2}\left \| \mathcal{C}-\mathcal{C}^{k}    \right \| _{F}^{2}, \end{cases} 
\end{equation}
where $ \rho _{1} $, $ \rho _{2} $ and $ \rho _{3} $ are given parameters, and $ \rho _{1},\rho _{2},\rho _{3> 0} $.	
It can be seen that all subproblems are strongly convex optimization problems, the existence and uniqueness of solutions are guaranteed, and all of them have explicit solutions. The details are as follows:\\
$ \mathcal{X} $-subproblem:
\begin{equation}\label{eq:3}
\mathcal{X}^{k+1} =\left ( \rho _{1}\mathcal{X}^{k}+\mathcal{C}^{k} \ast _{v}\left ( \mathcal{Y}^{k}   \right )^{\ast }\right )\left ( \mathcal{Y}^{k}\ast _{v} \left ( \mathcal{Y}^{k}   \right )^{\ast } +\rho _{1}I\right ) ^{\dagger },   
\end{equation}	
$ \mathcal{Y} $-subproblem:
\begin{equation}\label{eq:4}
\mathcal{Y}^{k+1} =\left ( \left ( \mathcal{X}^{k+1}   \right )^{\ast } \ast _{v} \mathcal{X}^{k+1}+\rho _{2}I\right ) ^{\dagger } \left ( \left ( \mathcal{X}^{k+1} \right )^{\ast } \ast _{v} \mathcal{C}^{k}+\rho _{2} \mathcal{Y}^{k}      \right ) , 
\end{equation}
$ \mathcal{C} $-subproblem:
$$
\mathcal{C}^{k+1}=\underset{\mathcal{C}}{argmin }f\left ( \mathcal{X}^{k+1},\mathcal{Y}^{k+1},\mathcal{C}  \right )+\frac{\rho _{3} }{2}\left \| \mathcal{C}-\mathcal{C}^{k}    \right \| _{F}^{2} .
$$
The above problem is equivalent to the following equality constraint problem:	
\begin{equation}\label{eq:5}
\begin{aligned}
\mathcal{C}^{k+1}\!=&\underset{\mathcal{C},\mathcal{Q}_{1},\mathcal{Q}_{2}   }{argmin} \frac{1}{2} \left \| \mathcal{X}^{k+1} \!\ast _{v}\mathcal{Y}^{k+1}\!-\!\mathcal{C}    \right \| _{F}^{2}+\!\alpha _{1}\left \| \mathcal{Q}_{1}   \right \|_{l_{1} } +\!\alpha _{2}\left \| \mathcal{Q}_{2}   \right \|_{l_{1} }+\!\frac{\rho _{3} }{2}\left \| \mathcal{C}\!-\!\mathcal{C}^{k}    \right \| _{F}^{2}+\!\Phi \left ( \mathcal{C}  \right )\\
&s.t. \mathcal{Q} _{1}= \mathcal{D} _{1}\ast _{v}\mathcal{C} , \mathcal{Q} _{2}=\mathcal{C}\ast _{v} \mathcal{D} _{2} .
\end{aligned}      
\end{equation}
Let $ \mathcal{S} $ and $ \mathcal{T} $ be Lagrangian multipliers of (\ref{eq:5}), then
\begin{equation}\label{eq:6}
\mathcal{Q} _{1}= \underset{\mathcal{Q} _{1}}{argmin} \alpha _{1}\left \| \mathcal{Q} _{1} \right \|_{l_{1} }+\frac{\beta }{2}\left \| \mathcal{Q} _{1}-\mathcal{D} _{1}\ast _{v}\mathcal{C}+\frac{\mathcal{S} }{\beta }    \right \| _{F}^{2}    , 
\end{equation}	
\vspace{-0.5cm}
\begin{equation}\label{eq:7}
\mathcal{Q} _{2}= \underset{\mathcal{Q} _{2}}{argmin} \alpha _{2}\left \| \mathcal{Q} _{2} \right \|_{l_{1} }+\frac{\mu }{2}\left \| \mathcal{Q} _{2}-\mathcal{C}\ast _{v}\mathcal{D} _{2}+\frac{\mathcal{T} }{\mu }    \right \| _{F}^{2}      ,
\end{equation}		
\begin{equation}\label{eq:8}
\begin{aligned}
\mathcal{C}^{k+1}=& \underset{\mathcal{C}}{argmin}\frac{1}{2} \left \| \mathcal{X}^{k+1}\ast _{v}\mathcal{Y}^{k+1} -\mathcal{C}      \right \|_{F}^{2}+\frac{\beta }{2} \left \| \mathcal{Q}_{1}-\mathcal{D}_{1} \ast _{v} \mathcal{C}+\frac{\mathcal{S} }{\beta }    \right \| \\
&+\frac{\mu }{2}\left \| \mathcal{Q}_{2}-\mathcal{C}\ast _{v}\mathcal{D} _{2}+\frac{\mathcal{T} }{\mu }\right \|_{F}^{2}+\frac{\rho _{3} }{2}\left \| \mathcal{C}-\mathcal{C}^{k} \right \|  +\Phi \left ( \mathcal{C}  \right )  ,  
\end{aligned}      
\end{equation}
\begin{equation}\label{eq:9}
\mathcal{S} =\mathcal{S}+\beta \left ( \mathcal{Q}_{1}-\mathcal{D}_{1} \ast _{v} \mathcal{C}     \right )       , 
\end{equation}	
\begin{equation}\label{eq:10}
\mathcal{T} =\mathcal{T}+\mu  \left ( \mathcal{Q}_{2}-\mathcal{C} \ast _{v} \mathcal{D}_{2}     \right )       . 
\end{equation}	
According to the soft threshold, equations (\ref{eq:6}) and (\ref{eq:7}) have the following unique solutions:	
\begin{equation}\label{eq:11}
\mathcal{Q}_{1}\left ( i,j,k \right )=T_{\frac{\alpha _{1} }{\beta } } \left ( \mathcal{D}_{1}\ast _{v}\mathcal{C} -\frac{\mathcal{S} }{\beta }     \right )\left ( i,j,k \right )     ,
\end{equation}	
\begin{equation}\label{eq:12}
\mathcal{Q}_{2}\left ( i,j,k \right )=T_{\frac{\alpha _{2} }{\mu } } \left (\mathcal{C}\ast _{v}\mathcal{D}_{2} -\frac{\mathcal{T} }{\mu }     \right )\left ( i,j,k \right )     ,
\end{equation}
where
$$
T_{\eta }\left ( x \right ) =\begin{cases} \left ( \left | x \right |-\eta   \right )sign\left ( x \right ),\left | x \right |> \eta    \\ 0,\left | x \right |\le \eta ,   \end{cases}  
$$
where $ x\in \mathbb{R} $, and formula (\ref{eq:8}) is equivalent to
\begin{equation}
\begin{aligned}
\bar{C}^{k+1}=&\underset{\bar{C}}{argmin}\frac{1}{2v}\left \| \bar{X}^{k+1}\bar{Y}^{k+1}-\bar{C}  \right \|_{F}^{2} +\frac{\beta }{2v}\left \| \bar{Q }_{1}-\bar{D}_{1}\bar{C}+\frac{\bar{S}}{\beta }   \right \|_{F}^{2}\\
&+\frac{\mu }{2v}\left \| \bar{Q }_{2}-\bar{C}\bar{D}_{2}+\frac{\bar{T}}{\mu }  \right \| +\frac{\rho _{3} }{2v}\left \| \bar{C}-\bar{C}^{k}   \right \| _{F}^{2} +\Phi \left (  \bar{C}\right )   .             
\end{aligned}
\end{equation}
The only solution of the above formula is actually the solution of the following matrix equation:
\begin{equation}\label{eq:13}
\bar{C}+\beta \bar{D}_{1}^{\ast } \bar{D}_{1}\bar{C}+\mu \bar{C}\bar{D}_{2}\bar{D}_{2}^{\ast }+\rho _{3}\bar{C}=R^{k}  ,   
\end{equation}
where
$$
R^{k}=\bar{X}^{k+1}\bar{Y}^{k+1}+\bar{D}_{1}^{\ast }\left ( \beta\bar{Q}_{1}+\bar{S}   \right )+\left ( \mu\bar{Q}_{2}+\bar{T}  \right )\bar{D}_{2}^{\ast }+\rho _{3}\bar{C}^{k}    .     
$$
The matrices in formula (\ref{eq:13}) all have block diagonal structures. From the definitions of $ \mathcal{D}_{1} $ and $ \mathcal{D}_{2} $, it can be seen that the diagonal blocks of $ \bar{D}_{1} $ and $ \bar{D}_{2} $ are the same, so this formula is equivalent to $ v $ matrix equations with smaller sizes, namely
\begin{equation}\label{eq:14}
\bar{C}_{l}+\beta H_{m}\bar{C}_{l}+\mu \bar{C}_{l} H_{n}+\rho _{3}\bar{C}_{l}=R_{l}^{k},l=1,2,\cdots ,v,       
\end{equation}
where $ \bar{C}_{l} $ and $ R_{l}^{k} $ represent the $ l $-th diagonal blocks of $ \bar{C}  $ and $ R^{k} $, respectively, and
$$
H_{m}:=L_{m}^{T}L_{m}=\begin{pmatrix}
1 & -1 &  &  & \\
-1 & 2 & -1 &  & \\
& \ddots  & \ddots  & \ddots  & \\
&  & -1 & 2 & -1\\
&  &  & -1 & 1
\end{pmatrix}\in \mathbb{R}^{m\times m},H_{m}^{0}:=I_{m} .    
$$
\begin{figure*}
	\centering
	\includegraphics[width=0.8\linewidth]{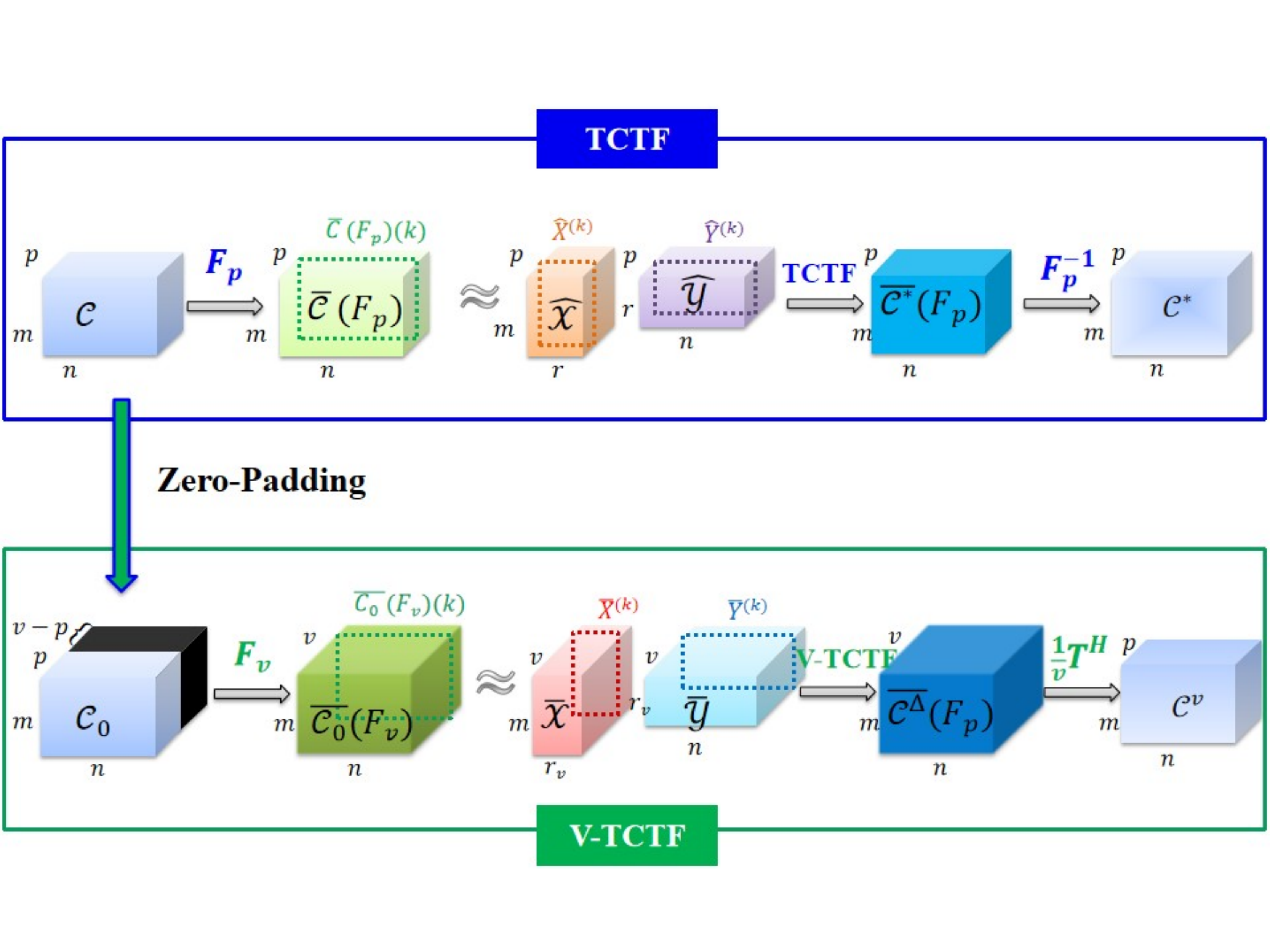}\\
	\caption{The difference between TCTF and V-TCTF}\label{fig:1}
\end{figure*}
It is easy to verify that $ H_{m} $ has the following orthogonal diagonalization form:
$$
H_{m}= K_{m}\Lambda _{m}K_{m}^{T},
$$
$$
K_{m}=\sqrt{\frac{2}{m} }\left [ \sqrt{\left ( 1+\delta _{j,1}  \right )^{-1}  } cos\left ( \frac{\pi \left ( 2i-1 \right )\left ( j-1 \right )  }{2m}  \right )  \right ]_{i,j=1}^{m},
$$
$$
\delta _{j,1}=\begin{cases} 1,j=1\\ 0,otherwise\end{cases},\\
\Lambda _{m}=4\times diag\left ( sin^{2}\left ( \frac{\left ( i-1 \right )\pi  }{2m}  \right )   \right )_{i=1}^{m}   .
$$
That is, $ \Lambda _{m} $ is a nonnegative diagonal matrix of $ m\times m $, so equation (\ref{eq:14}) is equivalent to
\begin{equation}\label{eq:15}
\bar{C}_{l}+\beta K_{m}\Lambda _{m}K_{m}^{T}  \bar{C}_{l}+\mu \bar{C}_{l}K_{n}\Lambda _{n}K_{n}^{T}  +\rho _{3}\bar{C}_{l}=R_{l}^{k},l=1,2,\cdots ,v.    
\end{equation}
Multiply $ K_{m}^{T} $ by left and $ K_{n} $ by right on both sides of equation (\ref{eq:15}) respectively to obtain
\begin{equation}\label{eq:16}
\left ( 1+\rho _{3}  \right ) \hat{C}_{l}+\beta \Lambda _{m}\hat{C}_{l}+\hat{C}_{l}\mu \Lambda _{n}=\hat{R}_{l}^{k} ,l=1,2,\cdots ,v,      
\end{equation}
where
$$
\hat{C}_{l}=K_{m}^{T}\bar{C}_{l}K_{n} ,\hat{R}_{l}^{k}=K_{m}^{T}R_{l}^{k} K_{n}   ,    
$$
and
\begin{equation}
\hat{C}_{l}\left ( i,j \right ) =\frac{\hat{R}_{l}^{k}\left ( i,j \right ) }{1+\rho _{3} +\beta \Lambda _{m}\left ( i \right )+\mu \Lambda _{n}\left ( j \right )    },\left ( i,j \right ) \in \Xi \left ( m \right )\times \Xi \left ( n \right ).
\end{equation}
By using
$$
\bar{C}_{l}=K_{m}\hat{C}_{l}K_{n}^{T},
$$
then we obtain
\begin{equation}\label{eq:17}
\mathcal{C}^{k+1}\left ( i,j,k \right )=\begin{cases} ifft_{v2p}\left ( \bar{\mathcal{C} } ,\left [  \right ],3   \right ),\left ( i,j,k \right )\notin \Omega    \\ \mathcal{G}\left ( i,j,k \right ),\left ( i,j,k \right )\in \Omega .   \end{cases}   
\end{equation}
\par For clarity, the algorithmic framework is presented in \cref{alg:1}.

\par When $ \alpha _{1} =\alpha _{2}=0 $, the model is simplified as:
\begin{equation}\label{opt:ZTCTF}
\begin{aligned}
\mathop {\min }\limits_{ {\cal X},  {\cal Y},{\C}} \quad & \frac{1}{2}  \left\| {\cal X}*_v {\cal Y} -  \C \right\|_F^2 \\
\mbox{s.t.} \quad  & \mathcal{C}_{\Omega }=\mathcal{G}_{\Omega }\\
&  \X \in \mathbb{R}^{m\times q\times p},   \Y \in \mathbb{R}^{q\times n \times p},\mathcal{C}\in \mathbb{R}^{m\times n\times p}.
\end{aligned}
\end{equation}
\par Using the same computational scheme proposed in \cite{ZLLZ18} to solve the Zero-Padding TCTF (\ref{opt:ZTCTF}) in the variable Fourier domain, we term it as V-TCTF. In summary, the difference between TCTF and our proposed V-TCTF is illustrated in Fig. \ref{fig:1}.
\begin{algorithm}
	\caption{The VTCTF-TV algorithm to solve (\ref{eq:2})}
	\label{alg:1}
	\begin{algorithmic}
		\STATE{\textbf{Input:} The tensor data $\mathcal{G}\in \mathbb{R}^{m\times n \times p}$, the observed set $\Omega$, the parameter $v=2p-1$, the initialized rank $r^0 \in \mathbb{R}^{v}$, and $\varepsilon=1e-5$, iteration $ N $=200.}
		\STATE{\textbf{Initialize:} ${\cal X}^{0}$, ${\cal Y}^{0}$;}
		\STATE{\textbf{While $ k\le N $ and not converged do}}
		\STATE{1. Fix $ \mathcal{Y}^{k} $ and $ \mathcal{C}^{k} $ to update $ \mathcal{X}^{k+1} $ via (\ref{eq:3}).}
		\STATE{2. Fix $ \mathcal{X}^{k+1} $ and $ \mathcal{C}^{k} $ to update $ \mathcal{Y}^{k+1} $ via (\ref{eq:4}).}
		\STATE{3. Fix $ \mathcal{X}^{k+1} $ and $ \mathcal{Y}^{k+1} $ to update $ \mathcal{C}^{k+1} $ via (\ref{eq:17}).}
		\STATE{4. Check the  termination criterion:$ \frac{\left \| \mathcal{C}^{k+1}-\mathcal{C}^{k}     \right \|_{F}^{2}  }{\left \| \mathcal{C}^{k+1}   \right \|_{F}^{2}  }\le \varepsilon $.}
		\STATE{\textbf{end while}}
		\STATE{\textbf{\textbf{Output:} $ \mathcal{C}^{k+1}$.}}
	\end{algorithmic}
\end{algorithm}

\par Obviously, the objective function in model (\ref{eq:1}) is proper lower semi-continuous, and only contains a simple projection constraint condition. According to the algorithm design, the generated point sequence $ \left \{ \left ( \mathcal{X}^{k},\mathcal{Y} ^{k} ,\mathcal{C}^{k}    \right )  \right \}  $ meets the conditions H1, H2 and H3 in Section 2.3 of reference\cite{attouch2013convergence}, so it can be proved that it converges to the critical point of model (\ref{eq:1}). The concrete convergence analysis is given in \cref{sec:Convergence analysis}.

\section{Numerical Experiments}\label{sec:Numerical}
In this section, Algorithm \ref{alg:1} is used to restore color images and multispectral images to evaluate its performance, and compared with the following four data completion methods, namely TCTF\cite{ZLLZ18}, Tmac\cite{xu2013parallel}, TCTFTVT\cite{lin2020tensor} and MTRTC\cite{yu2020multi}. In order to quantitatively evaluate the image quality restored by each method, PSNR\cite{chen2015fractional}, SSIM\cite{liu2014generalized} and CPU time are used as numerical indicators, and PSNR and SSIM are defined as follows:
\begin{equation*}
\begin{split}
&PSNR=10\cdot log_{10} \frac{mnp\left \| \mathcal{C}_{true}   \right \|_{\infty }^{2}  }{\left \| \mathcal{C}-\mathcal{C}_{true}    \right \| _{F}^{2} },\\
&SSIM=\frac{\left ( 2\mu _{\mathcal{C} }\mu _{\mathcal{C}_{true}  }   \right )\left ( 2\sigma _{\mathcal{C}\mathcal{C}_{true}   }+c_{2}   \right )  }{\left ( \mu _{\mathcal{C}  }^{2}\mu _{\mathcal{C}_{true}  }^{2}+c_{1} \right ) \left ( \sigma _{\mathcal{C} }^{2} +\sigma _{\mathcal{C}_{true}  }^{2} +c_{2}  \right ) }  ,
\end{split}
\end{equation*}
where $ \mathcal{C}_{true} $ is the real tensor, $ \mathcal{C} $ is the recovery tensor, $ \mu _{\mathcal{C} } $ and $ \mu _{\mathcal{C}_{true}  } $ are the mean values of images $ \mathcal{C} $ and $ \mathcal{C}_{true} $ , $ \sigma _\mathcal{C} $ and $ \sigma _{\mathcal{C}_{true}  } $ are the standard deviations of images $ \mathcal{C} $ and $ \mathcal{C}_{true} $ , $ \sigma _{\mathcal{C}\mathcal{C}_{true}   } $ is the covariance of images $ \mathcal{C} $ and $ \mathcal{C}_{true} $ , $ c_{1} $ and $ c_{2}> 0 $ are constants. The higher the $ PSNR $ and $ SSIM $ values, the better the image quality.
\par The parameters of the VTCTF-TV method proposed in this paper are set as follows:
$$
\alpha _{1}=\alpha _{2}=\beta=\mu =\varepsilon=1e-5,\rho _{1}=\rho _{2}=\rho _{3}=5e-6,N=200 .
$$
Set the initial rank to $ r^{0}=\left [ 30,\cdots ,30 \right ] \in \mathbb{R}^{v} $, the experiment is implemented on Matlab R2016b under Windows 10, equipped with a 3.00GHz CPU and a PC with 8GB of memory.

\subsection{Gray video}\label{subsec:N1}
In order to find the optimal value of the introduced additional parameter $ v $, we evaluate our method on the widely used YUV Video Sequences\footnote{ http://trace.eas.asu.edu/yuv/}. Each sequence contains at least 150 frames. In the experiments, we first test V-TCTF on the Hall Monitor video. The frame size of this video is 144$ \times  $176 pixels. Due to the computational limitation, we only use the first 20 frames of the sequences, thus p = 20 here. The sampling rate is set to 0.3, 0.5 and 0.7, and the initialized rank is set to $ r^{0}=\left [ 30,\cdots ,30 \right ] \in \mathbb{R}^{v} $. 
\begin{figure}[H]
	$$
	\begin{array}{cc}
	\includegraphics[width=0.47\textwidth,height=0.25\textwidth]{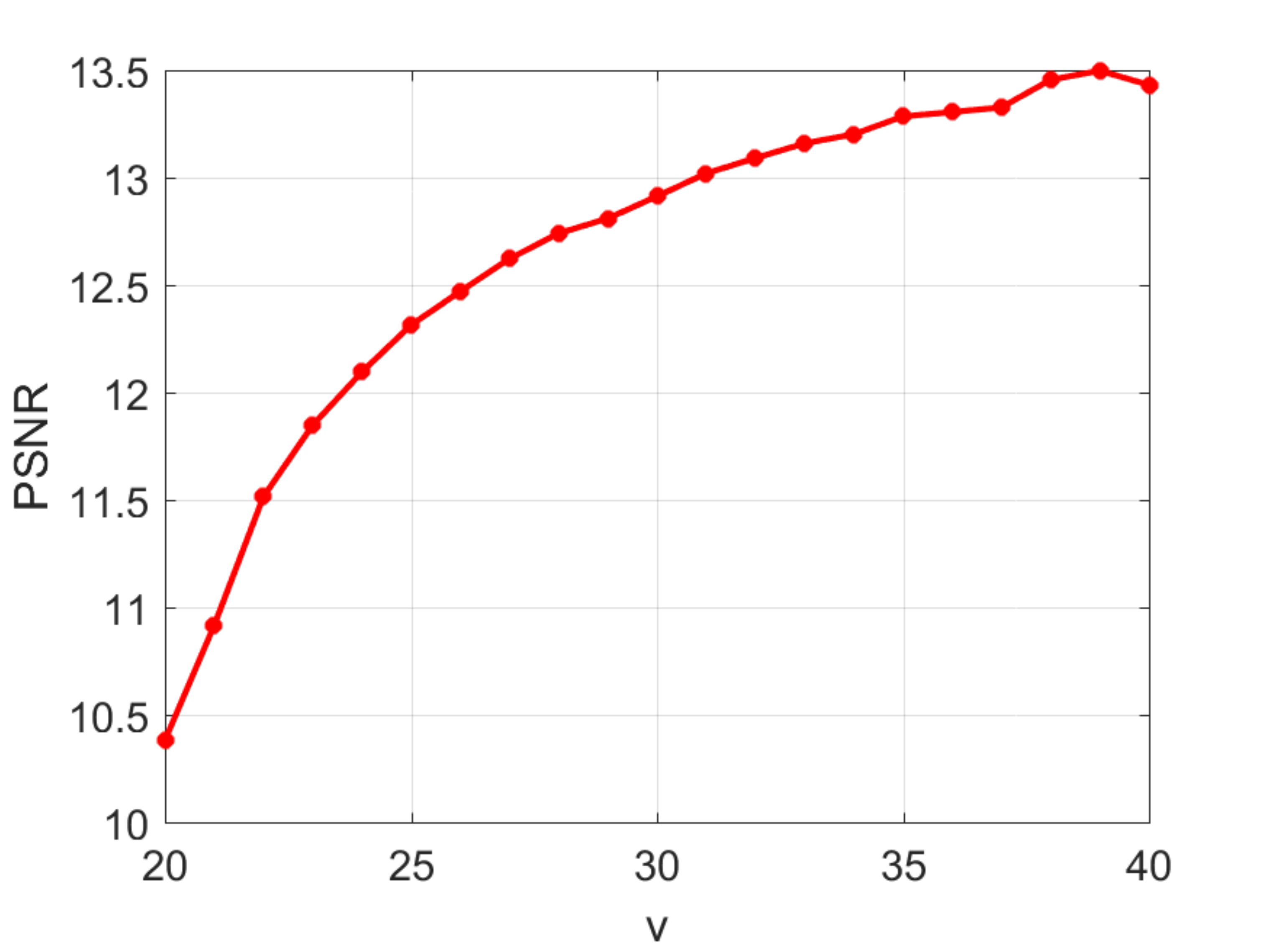}&
\includegraphics[width=0.47\textwidth,height=0.25\textwidth]{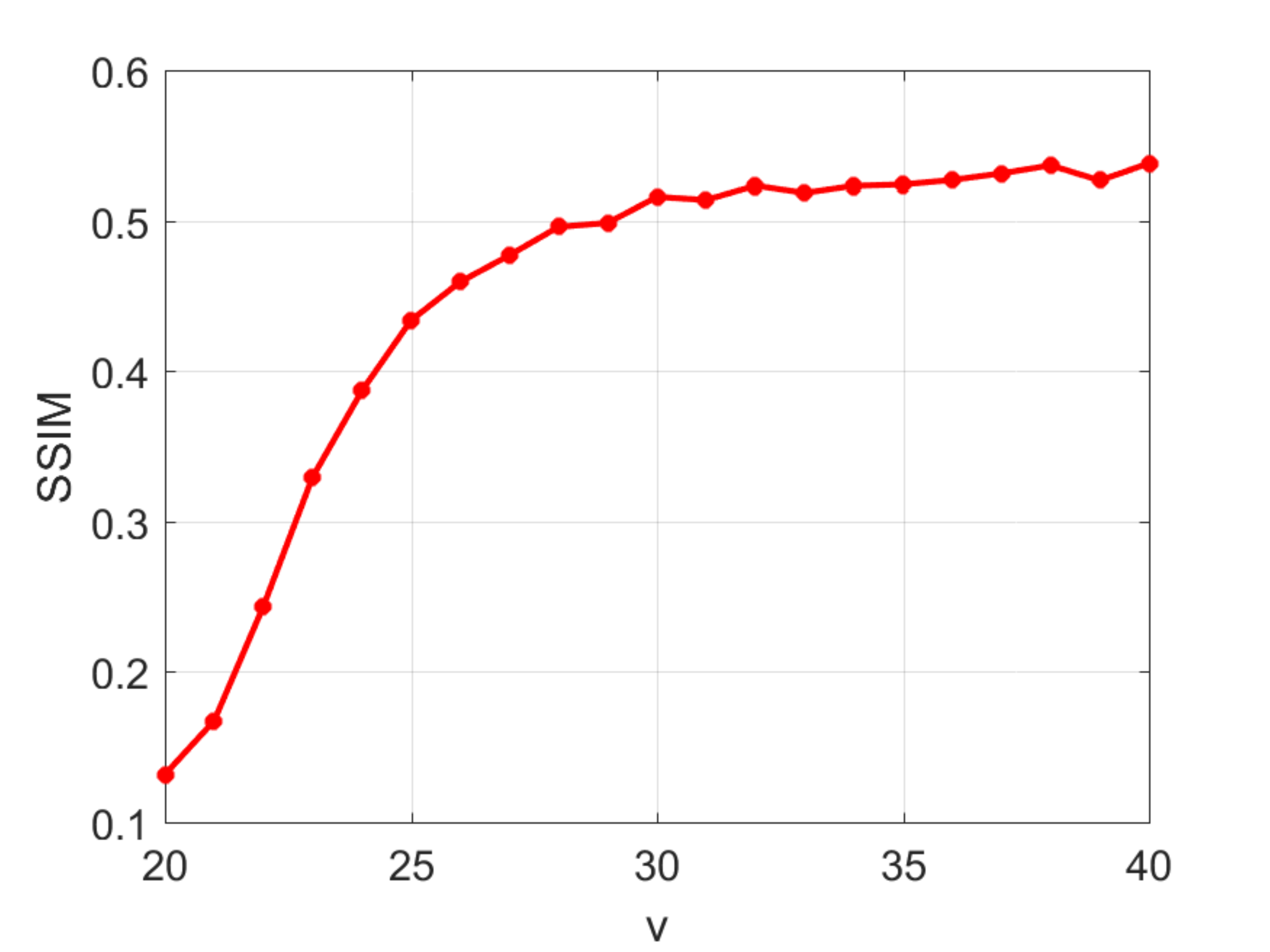}\\
\includegraphics[width=0.47\textwidth,height=0.25\textwidth]{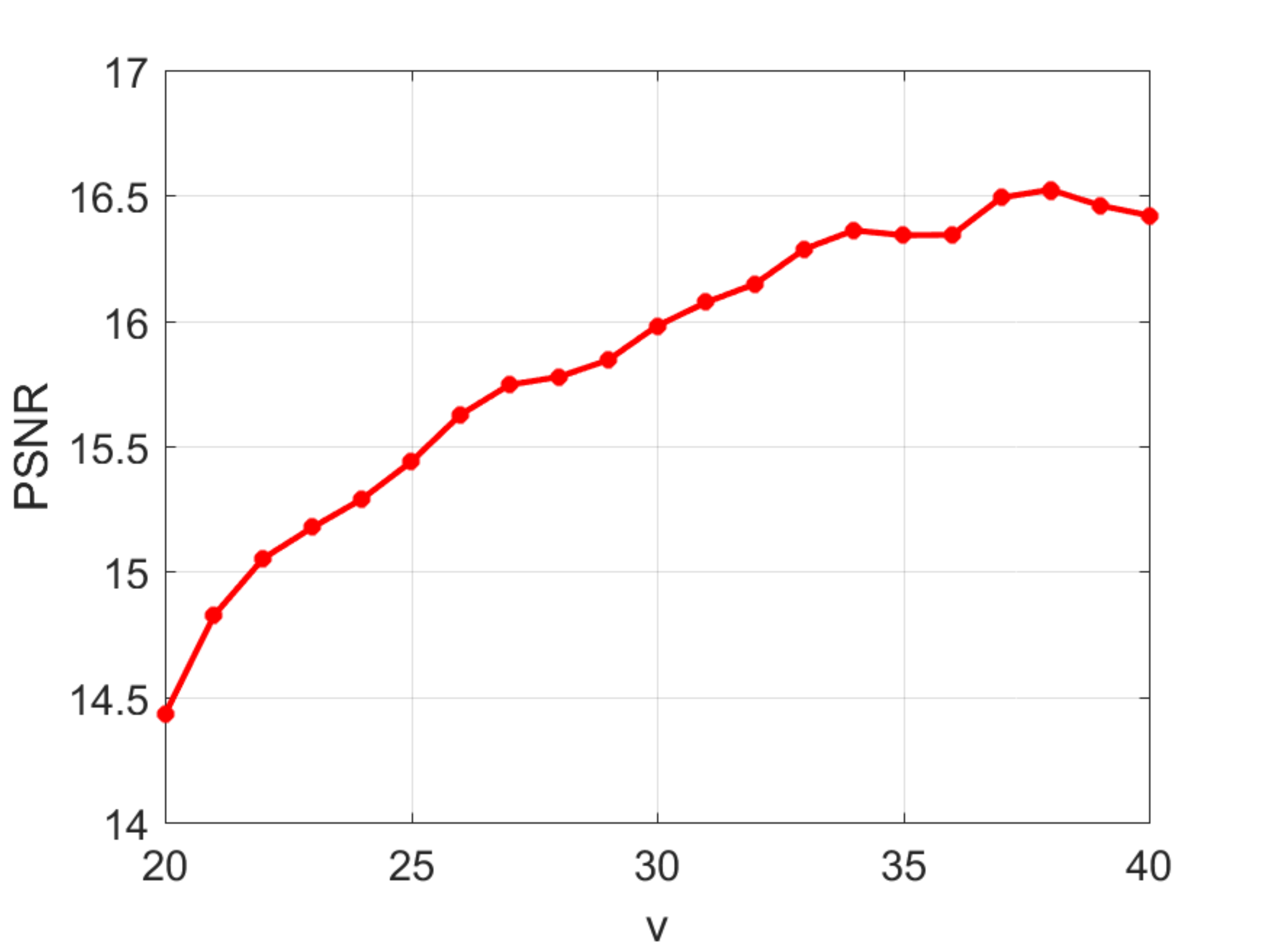}&
\includegraphics[width=0.47\textwidth,height=0.25\textwidth]{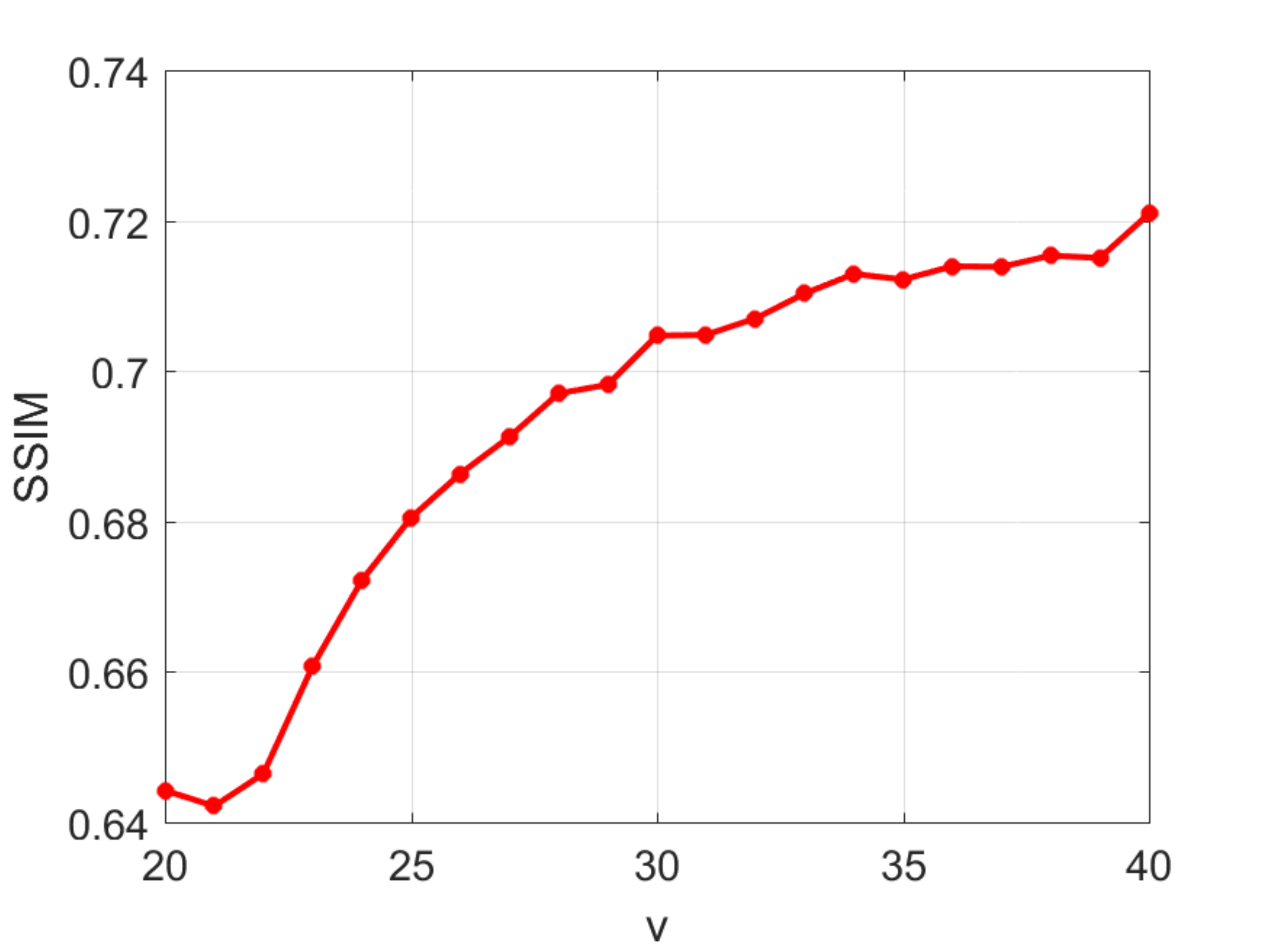}\\
\includegraphics[width=0.47\textwidth,height=0.25\textwidth]{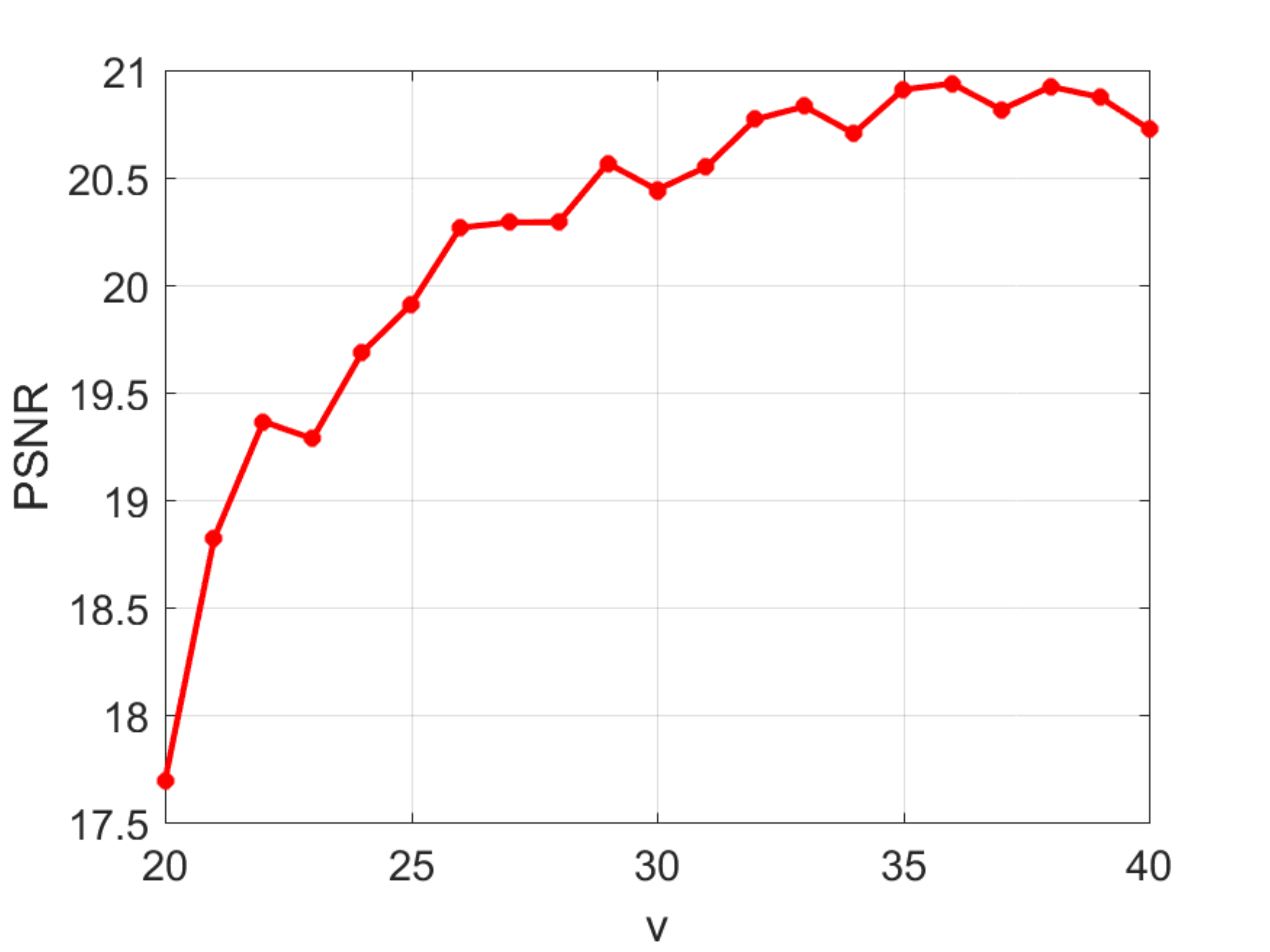}&
\includegraphics[width=0.47\textwidth,height=0.25\textwidth]{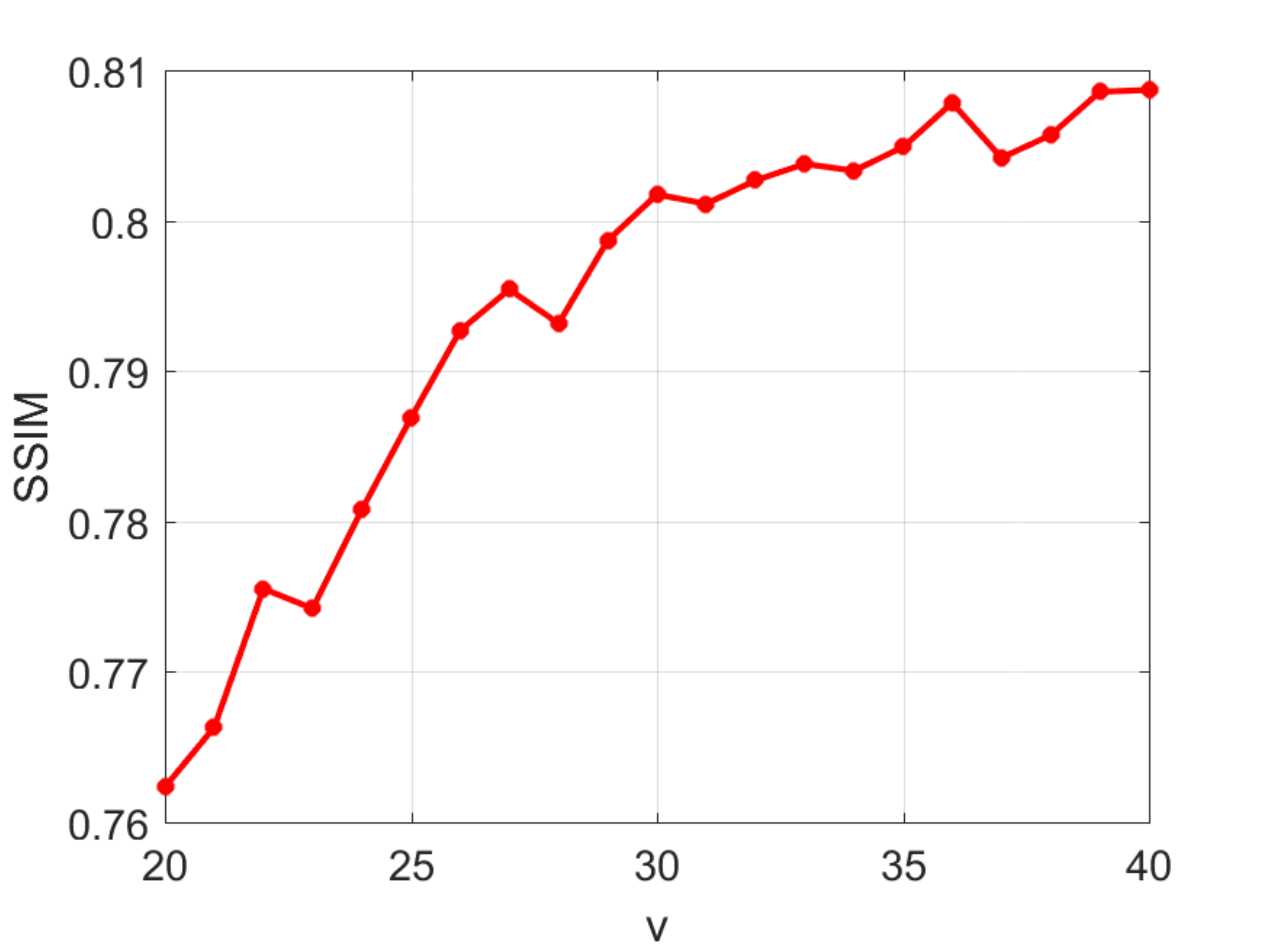}\\
	\end{array}
	$$
	\caption{The PSNR and SSIM value v.s. different v for V-TCTF method on the hall video, from top to bottom, the sampling rate is 0.3, 0.5 and 0.7.}\label{fig:2}	
\end{figure} 
As shown in Fig. \ref{fig:2}, we display the PSNR and SSIM comparison of the testing video. We can see from the Fig. \ref{fig:2} that the resulting tensor completion can be performed much better when $ v=2p-1 $ or near $ 2p $.
\subsection{Color images}\label{subsec:N2}
In this subsection, we conduct experiments on six popular color images (Lena, Panda, Sailboat, Barbara, House and Pepper) that are widely used in the literature. All images are of size 256 $\times$ 256 $\times$ 3. We compared the performance of our method with five other image completion methods, including TCTF method, TMac, TCTFTVT method and MTRTC method.
\begin{figure}[H]
	\centering
	\includegraphics[width=1.0\linewidth,height=0.98\linewidth]{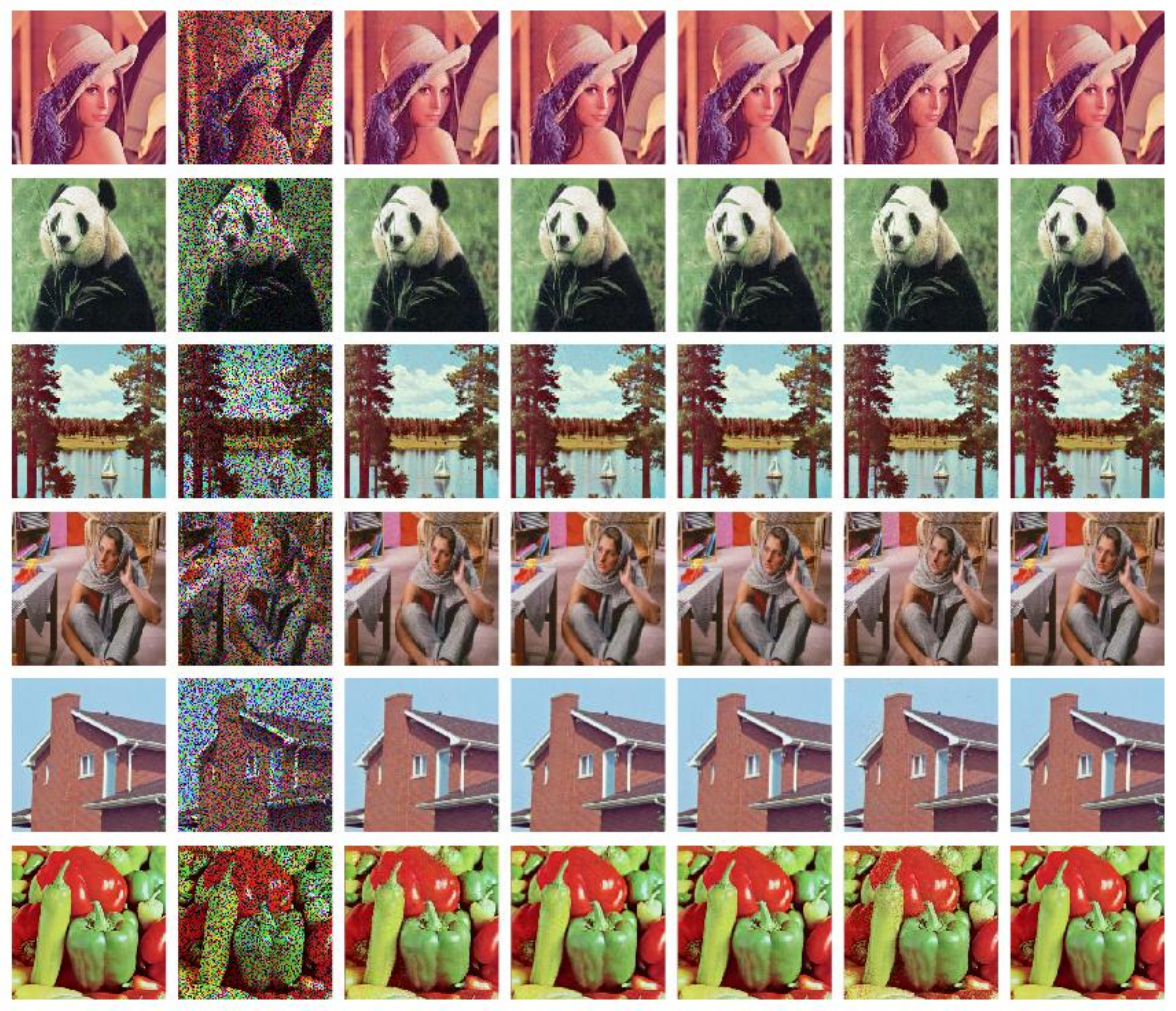}\\
	\caption{Restored results of color images with sampling rate 0.7. From top to bottom, Lena, Panda, Sailboat, Barbara, House and Pepper. From left to right: the original data, the observed data, the recovered
		results by TCTF, Tmac, TCTFTVT, MTRTC and VTCTF-TV, respectively.}\label{fig:3}
\end{figure}
\vspace{-0.9cm}
In Table \ref{tab:1}, by comparing the peak signal-to-noise ratio, structural similarity and operation time at different sampling rates of 0.6, 0.7 and 0.8, it is easy to see that the peak signal-to-noise ratio and structural similarity of VTCTF-TV are higher than those of TCTF, TMac, TCTFTVT and MTRTC.
\begin{table}[H]
	\centering
	\normalsize
	\begin{spacing}{1.5} 
		\resizebox{\textwidth}{!}{ 
			\begin{tabular}{ccccccccc}
				\toprule
				Color image & SR & TCTF & Tmac & TCTFTVT & MTRTC & VTCTF-TV \\
				\midrule
				\multirow{3}{*}{Lena} & 0.6 & 28.40/0.83/4.10	&	29.90/0.86/9.00	&	30.18/0.90/7.63	&	29.47/0.84/3.31	&	\textbf{30.60}/\textbf{0.93}/2.78\\
				& 0.7 & 29.78/0.86/3.68	&	30.79/0.88/8.73	&	31.52/0.92/7.83	&	30.77/0.86/3.10 &	\textbf{32.47}/\textbf{0.96}/3.85\\
				& 0.8 & 31.68/0.91/3.23	&	31.69/0.90/8.49	&	33.35/0.95/7.80	&	33.15/0.91/2.76	&	\textbf{34.87}/\textbf{0.97}/3.07\\
				\cline{2-7}
				\multirow{3}{*}{Panda} & 0.6 & 30.40/0.85/3.25	&	30.99/0.83/7.74	&	31.67/0.88/8.29	&	31.28/0.87/3.86	&	\textbf{31.88}/\textbf{0.90}/3.70\\
				& 0.7 & 31.83/0.88/2.97	&	31.92/0.86/8.84	&	33.00/0.91/7.89	&	33.14/0.91/2.83 &	\textbf{33.32}/\textbf{0.92}/3.10\\
				& 0.8 & 33.70/0.92/3.51	&	32.79/0.88/8.44	&	34.79/0.94/7.89	&	35.32/0.94/2.83	&	\textbf{35.65}/\textbf{0.95}/3.10\\
				\cline{2-7}
				\multirow{3}{*}{Sailboat} & 0.6 & 25.73/0.79/3.24	&	27.32/0.81/7.66	&	27.48/0.85/7.28	&	26.97/0.83/3.83	&	\textbf{27.76}/\textbf{0.90}/14.98\\
				& 0.7 & 27.16/0.84/3.75	&	28.21/0.84/9.68	&	28.86/0.89/7.22	&	28.68/0.87/3.45 &	\textbf{29.81}/\textbf{0.94}/4.10\\
				& 0.8 & 29.06/0.89/3.68	&	29.13/0.87/8.46	&	30.74/0.92/7.60	&	30.92/0.91/2.79	&	\textbf{32.17}/\textbf{0.96}/2.93\\
				\cline{2-7}
				\multirow{3}{*}{Barbara} & 0.6 & 27.85/0.84/3.42	&	29.41/0.86/8.09	&	\textbf{29.73}/0.89/7.45	&	28.16/0.86/3.41	&	29.70/\textbf{0.91}/14.81\\
				& 0.7 & 29.27/0.87/3.80	&	30.26/0.88/8.53	&	31.10/0.92/7.98	&	29.37/0.88/3.16 &	\textbf{32.01}/\textbf{0.94}/3.78\\
				& 0.8 & 31.13/0.91/3.38	&	31.14/0.90/8.59	&	32.96/0.95/8.04	&	31.57/0.91/2.63	&	\textbf{34.37}/\textbf{0.97}/3.02\\
				\cline{2-7}
				\multirow{3}{*}{House} & 0.6 & 30.52/0.86/2.89	&	33.09/0.89/7.84	&	33.18/0.91/8.12	&	28.98/0.80/3.12	&	\textbf{33.48}/\textbf{0.92}/4.63\\
				& 0.7 & 31.90/0.89/3.83	&	34.04/0.91/8.74	&	34.58/0.94/7.53	&	30.98/0.85/2.84 &	\textbf{35.56}/\textbf{0.95}/3.39\\
				& 0.8 & 33.73/0.93/3.62	&	34.99/0.92/9.31	&	36.39/0.96/7.66	&	33.32/0.90/2.65	&	\textbf{37.83}/\textbf{0.96}/2.92\\
				\cline{2-7}
				\multirow{3}{*}{Pepper} & 0.6 & 25.15/0.72/3.37	&	28.97/0.84/8.77	&	29.48/0.88/8.18	&	24.49/0.71/3.12	&	\textbf{30.66}/\textbf{0.93}/5.11\\
				& 0.7 & 26.57/0.77/3.59	&	29.92/0.87/8.31	&	30.79/0.91/8.26	&	25.72/0.75/2.83 &	\textbf{33.02}/\textbf{0.96}/3.90\\
				& 0.8 & 28.41/0.83/3.27	&	30.96/0.89/8.94	&	32.69/0.94/8.47	&	27.45/0.80/2.61	&	\textbf{35.30}/\textbf{0.97}/3.39\\
				\bottomrule
			\end{tabular}
		}
	\end{spacing}
	\caption{The PSNR/SSIM/CPU time of different methods on color images ``Lena", ``Panda", ``Sailboat", ``Barbara", ``House" and ``Pepper" with different sampling rates.}\label{tab:1}
\end{table}
\vspace{-0.8cm}
 Fig. \ref{fig:3} is a restored image of six color images in models TCTF, TMac, TCTFTVT, MTRTC and VTCTF-TV with the sampling rate set to 0.7. 
 It can be observed that the recovery effect of VTCTF-TV is better than that of TCTF, TMac, TCTFTVT and MTRTC.

\subsection{Multispectral images}\label{subsec:N3}
In this section, the performance of different methods on multi-spectral images is tested, and three multi-spectral images\footnote{https://www.cs.columbia.edu/CAVE/databases/multispectral/stuff/} are numerically tested, and their peak signal-to-noise ratio, structural similarity and CPU time are compared. The size of three multispectral images is 512 $ \times  $512 $ \times  $31.
\begin{figure}[H]
	\centering
	\includegraphics[width=1.0\linewidth,height=0.31\linewidth]{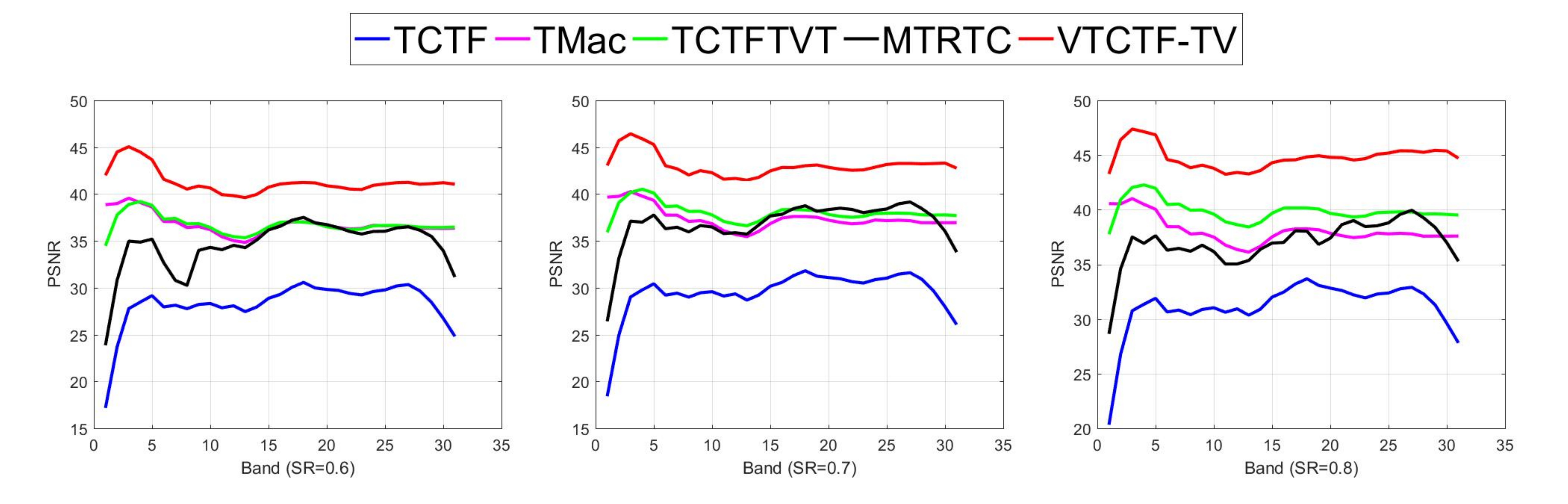}\\
	\caption{PSNR values of each band of the recovered MSI  Pompoms .}\label{fig:6}
\end{figure}
\begin{figure}[H]
	\centering
	\includegraphics[width=1.0\linewidth,height=0.31\linewidth]{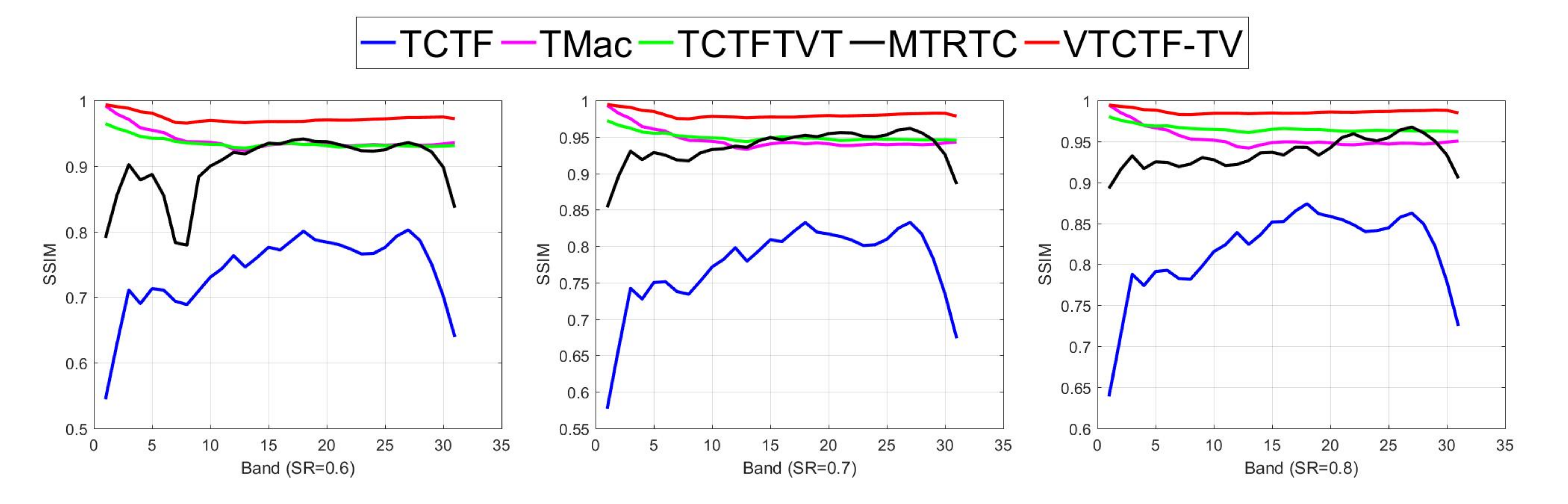}\\
	\caption{SSIM values of each band of the recovered MSI  Pompoms .}\label{fig:7}
\end{figure}
\vspace{-0.8cm}
\begin{figure}[H]
	\centering
	\includegraphics[width=1.001\linewidth]{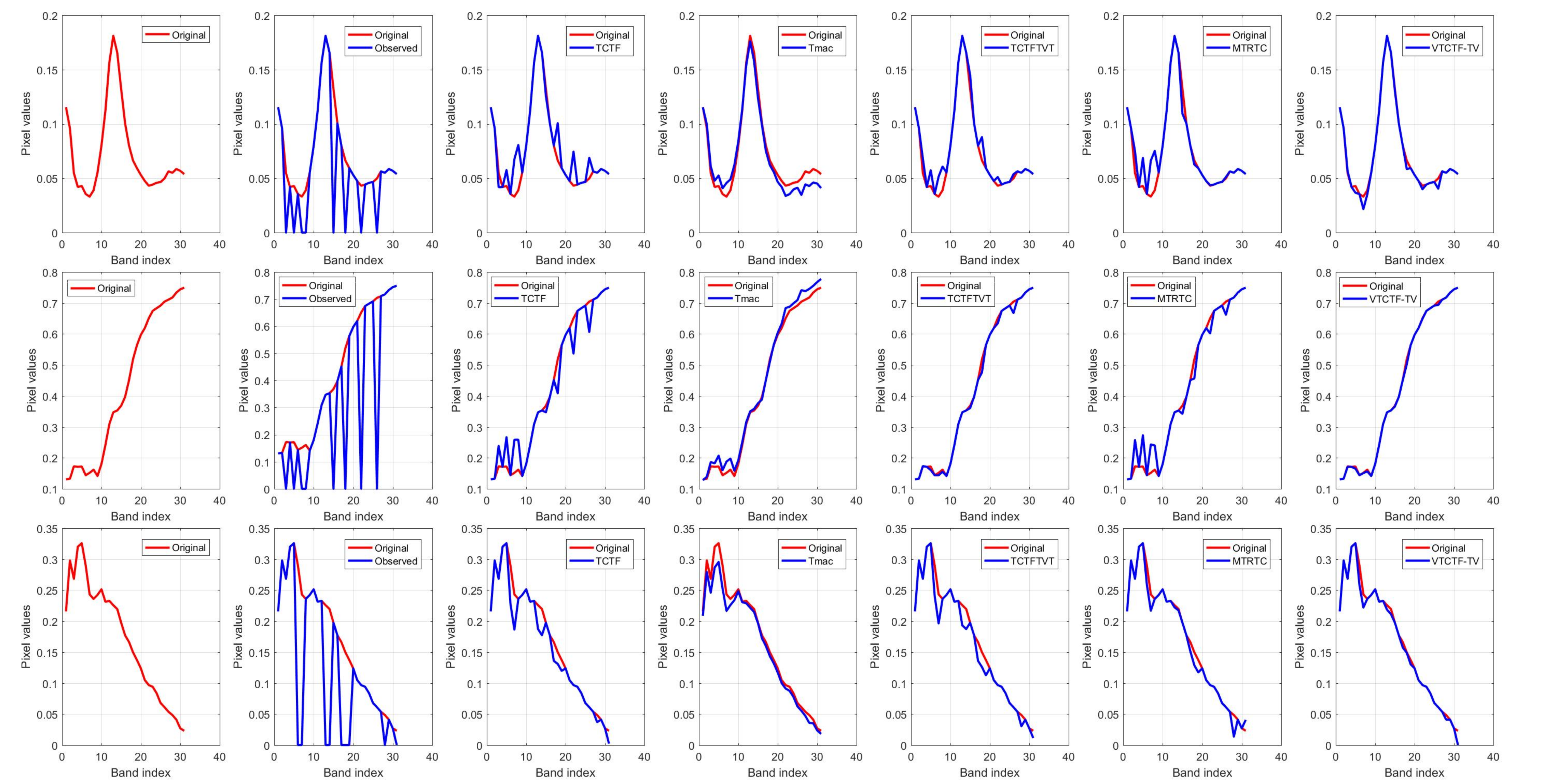}\\
	\caption{The one mode-3 tube of recovered MSIs by different methods with SR=0.7. From top to bottom: Pompoms ,  Flowers  and  Stuffed toys, respectively. From left to right: the original data, the observed data, the recovered results by TCTF, Tmac, TCTFTVT, MTRTC and VTCTF-TV, respectively.}\label{fig:5}
\end{figure}
\vspace{-0.7cm}
In Table \ref{tab:2}, by comparing the average peak signal-to-noise ratio, average structural similarity and CPU time at different sampling rates of 0.6, 0.7 and 0.8, it is easy to see that the peak signal-to-noise ratio and structural similarity of VTCTF-TV are higher than those of TCTF, TMac, TCTFTVT and MTRTC.  Fig. \ref{fig:6} and Fig. \ref{fig:7} give the PSNR and SSIM values comparison of each band of the MSI  Pompoms  recovered by all methods, respectively. The results on the other two datasets are similar.
\par In addition, Fig. \ref{fig:5} displays one mode-3 tube of the recovered MSIs by different methods with SR = 0.7, we clearly observe that 
VTCTF-TV yields the closest mode-3 tubes in all cases. 
\par
\begin{figure}[H]
	\centering
	\includegraphics[width=1.0\linewidth]{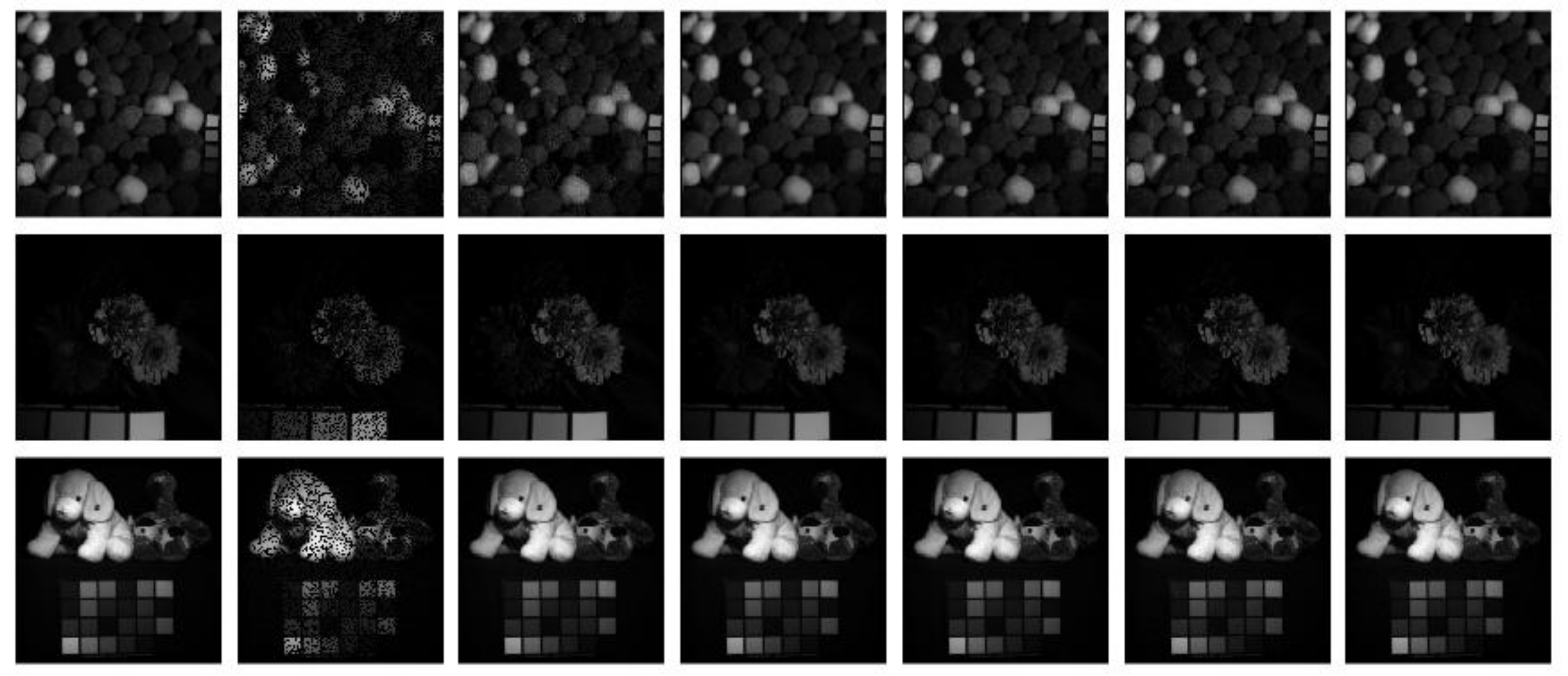}\\
	\caption{Restored results on MSIs with 0.7 observed entries. From top to bottom: The tenth band of MSI  Pompoms ,  Flowers  and  Stuffed toys. From left to right: the original data, the observed data, the recovered results by TCTF, Tmac, TCTFTVT, MTRTC and VTCTF-TV, respectively.}\label{fig:4}
\end{figure}

 Fig. \ref{fig:4} is the restored image of the 10th band of three multispectral images when the sampling rate of models TCTF, TMac, TCTFTVT, MTRTC and VTCTF-TV is set to 0.7. From left to right are the original picture, the observed image and the experimental results under TCTF, TMac, TCTFTVT, MTRTC and VTCTF-TV algorithms. From the visual comparison, it is clear that VTCTF-TV performs best in preserving multispectral image edges and details. In a word, VTCTF-TV has better recovery effect than TCTF, TMac, TCTFTVT and MTRTC.
\begin{table}[H]
	\centering
	\normalsize
	\begin{spacing}{1.5} 
		\resizebox{\textwidth}{!}{ 
			\begin{tabular}{ccccccccc}
				\toprule
				Multispectral image & SR & TCTF & Tmac & TCTFTVT & MTRTC & VTCTF-TV \\
				\midrule
				\multirow{3}{*}{ Pompoms } & 0.6 & 28.21/0.74/95.81	&	36.80/0.94/100.94	&	36.72/0.94/175.47	&	34.53/0.90/85.09	&	\textbf{41.31}/\textbf{0.97}/148.09\\
				& 0.7 & 29.47/0.77/123.47	&	37.42/0.95/122.71	&	38.04/0.95/178.49	&	36.84/0.94/75.43 &	\textbf{43.08}/\textbf{0.98}/96.18\\
				& 0.8 & 31.06/0.82/105.76	&	38.09/0.95/113.28	&	39.86/0.97/178.73	&	36.97/0.94/67.87	&	\textbf{44.80}/\textbf{0.99}/74.49\\
				\cline{2-7}
				
				\multirow{3}{*}{ Flowers } & 0.6 & 31.89/0.87/100.76	&	36.78/0.94/106.15	&	36.53/0.95/215.37	&	36.92/0.94/192.47	&	\textbf{39.17}/\textbf{0.97}/212.02\\
				& 0.7 & 32.78/0.89/139.82	&	37.41/0.95/151.71	&	37.91/0.97/247.35	&	37.54/0.95/140.43 &	\textbf{41.10}/\textbf{0.98}/108.00\\
				& 0.8 & 34.11/0.92/108.01	&	38.09/0.96/115.71	&	39.76/0.98/177.81	&	39.02/0.96/146.64	&	\textbf{43.50}/\textbf{0.99}/86.34\\
				\cline{2-7}
				
				\multirow{3}{*}{ Stuffed toys } & 0.6 & 39.73/0.97/107.05	&	40.87/0.97/121.59	&	40.83/\textbf{0.98}/156.97	&	37.10/0.92/107.25	&	\textbf{41.83}/0.96/116.93\\
				& 0.7 & 41.29/0.97/108.68	&	41.55/\textbf{0.98}/121.62	&	42.20/0.98/160.13	&	38.65/0.94/118.04 &	\textbf{43.53}/0.97/89.61\\
				& 0.8 & 34.97/0.97/116.43	&	42.28/0.98/126.03	&	44.06/0.99/162.50	&	43.96/0.98/118.74	&	\textbf{45.89}/\textbf{0.99}/73.23\\
				\bottomrule
			\end{tabular}
		}
	\end{spacing}
	\caption{The PSNR/SSIM/CPU time of different methods on multispectral images `` Pompoms ", `` Flowers " and `` Stuffed toys " with different sampling rates.}\label{tab:2}
\end{table}
\vspace{-0.8cm}
\vspace{-0.8cm}

\section{Conclusion}\label{sec:con}
In this paper, we propose variable T-product, and apply the Zero-Padding Discrete Fourier Transformation (ZDFT), to convert third order tensor problems to the variable Fourier domain, and a novel tensor completion model based on ZDFT and TV regularization is proposed. 
In the model, tensor decomposition based on variable T-product is used to describe the low rank of tensor, and TV regularization is used to describe the smoothness of image data.
Numerical experiments on color images and multispectral images show the effectiveness of the proposed algorithm.
\newpage
\appendix
\section{Proof of Proposition \ref{p3.1}}
\begin{proof} Let $\aa, \bb \in \RR^p$, and $\cc = \phi(\aa) \circ \phi(\bb)$.   For $k = 1, \cdots, p$, we have
\begin{eqnarray*}
	&& \phi^H(\cc)(k)\\
	& = & {1 \over v}\sum_{l=1}^v \bar \omega^{(l-1)(k-1)} \cc(l)\\
	& = & {1 \over v}\sum_{l=1}^v\sum_{i=1}^p\sum_{j=1}^p \bar \omega^{(l-1)(k-1)}\omega^{(i-1)(l-1)}\omega^{(j-1)(l-1)}\aa(i)\bb(j)\\
	& = & {1 \over v}\sum_{i=1}^p\sum_{j=1}^p \aa(i)\bb(j)\left[\sum_{l=1}^v\omega^{(l-1)(i+j-k-1)}\right].
\end{eqnarray*}
Note that
$$\sum_{l=1}^v\omega^{(l-1)(i+j-k-1)} = \left\{ \begin{aligned}v, & \ \ {\rm if\ } i+j-k-1= 0, \ {\rm mod}(v) \\0, & \ \ {\rm otherwise}. \end{aligned}
\right.$$
Thus,
$$\phi^H(\cc)(k) =  \sum \left\{ \aa(i)\bb(j) : \begin{array}{l}
i+j-k-1= 0\ {\rm mod}(v), \\
i, j =1, \cdots, p
\end{array}\right\}.$$
By (\ref{e2.1}), $\phi^H(\cc) = \aa \odot_v \bb$.
This proves (\ref{e3.3}).
\end{proof}

\section{Proof of Proposition \ref{truncated-T}}
\begin{proof} Denote $F_v = [T ~N]$, where $N\in \CC^{(v-p)\times v}$ consists of the last $v-p$ columns of $F_v$. Note that $F_v^{-1} = \frac{1}{v}F^H$. Thus,
\begin{equation}\label{conj}
F_v^{-1}([1:p],:) = \frac{1}{v}T^H.
\end{equation}
Let $\C : = \A *_v \B$ and ${\cal D}: = \A_0 * \B_0$. It follows from Proposition \ref{p3.1} that for any $i\in \{1,\cdots, m\}$ and any $j\in \{1,\cdots, n\}$,
\begin{eqnarray*}
	\cc_{ij} &=& \frac{1}{v} \phi^H \left[\sum\limits_{l=1}^q \phi(\aa_{il})\circ \phi(\bb_{lj})\right] \nonumber\\
	& = & \frac{1}{v} T^H \left[\sum\limits_{l=1}^q T \aa_{il} \circ T \bb_{lj} \right] \nonumber\\
	& = & \frac{1}{v} T^H \left[\sum\limits_{l=1}^q F_v  \tilde{\aa}_{il} \circ F_v \tilde{\bb}_{lj} \right] \nonumber\\
	& = & F_v^{-1}([1:p],:) \left[\sum\limits_{l=1}^q F_v  (\aa_0)_{il} \circ F_v (\bb_0)_{lj} \right] \nonumber\\
	& = & \left(F_v^{-1} \left[\sum\limits_{l=1}^q F_v  (\aa_0)_{il} \circ F_v (\bb_0)_{lj} \right]\right)(1:p) \nonumber\\
	& = & {\bf d}_{ij} (1:p), \nonumber
\end{eqnarray*}
where the third equality is from \eqref{zp} and the fourth one from \eqref{conj}. Thus, ${\cal C} = {\cal D}(:,:,[1:p])$. This completes the proof. 
\end{proof}

\section{Proof of Theorem \ref{t10}}
\begin{proof}
The equivalence in the first part can be derived by using
$$\bar{\mathcal{C}}(T) = \bar{{\cal X}} *_H \bar{{\cal Y}} \Longleftrightarrow \bar{C}(T)^{(k)} =\bar{ X}^{(k)}\bar{Y}^{(k)}, k=1,\cdots, v$$
from the definitional expression \eqref{e2.24}, and by the definition of the variable tubal rank as stated in Definition \ref{v-tubal-rank}. To show the ``furthermore" part, let $\A_0$ and $\B_0$ be the zero-padding counterparts of $\A$ and $\B$, respectively, and ${\cal D}$ be the T-product of $\A_0$ and $\B_0$, i.e., ${\cal D} := \A_0 * \B_0$. It is known from Proposition \ref{truncated-T} that Rank$_v(\C) \leq {\rm Rank}_t{\cal D}$. Additionally, by Lemma 2 in \cite{ZLLZ18}, we have ${\rm Rank}_t{\cal D} \leq \min\{ {\rm Rank}_t (\A_0), {\rm Rank}_t(\B_0)\}\leq r$. Thus, Rank$_v(\C) \leq r$. This completes the whole proof. 
\end{proof}

\section{Convergence analysis for \cref{alg:1} }\label{sec:Convergence analysis}
\par In this section, we will prove the global convergence of the proposed algorithm.  
First and foremost, for positive integers $ I_1 $,$ I_2 $,$ I_3 $, define bijections $ \mathbb{V}_{[ I_1,I_2,I_3] }:\mathcal{A}\in\mathbb{R}^{I_1\times I_2\times I_3}\mapsto\mathbb{V}_{[ I_1,I_2,I_3] }(\mathcal{A})\in\mathbb{R}^{I_1 I_2 I_3\times 1} $ and $ \mathbb{I}_{[ I_1,I_2,I_3] }:\left( i_1,i_2,i_3\right)\in \Xi\left( I_1\right)\times\Xi\left( I_2\right)\times\Xi\left( I_3\right)\mapsto\mathbb{I}_{\left[ I_1,I_2,I_3\right] }\left( i_1,i_2,i_3\right) \in\Xi\left( I_1 I_2 I_3\right) $ by:
$$
 \mathbb{V}_{\left[ I_1,I_2,I_3\right] }\left( \mathcal{A}\right): = \mathcal{A}\left( :\right) ,\\
 $$
 $$
 \mathcal{A}\left( i_1,i_2,i_3\right) = \left[ \mathcal{A}\left( :\right) \right] \left( \mathbb{I}_{\left[ I_1,I_2,I_3\right] }\left( i_1,i_2,i_3\right) \right) \mbox{holding for each } \mathcal{A}\in\mathbb{R}^{I_1\times I_2\times I_3}.
$$
\par Denote $ d_1=m q p$ , $d_2=q n p$ ,$ d_3=m n p $, $ d=d_1+d_2+d_3 $, and denote $\mathbb{V}_{\left[ m,q,p\right] } $, $ \mathbb{V}_{\left[ q,n,p\right] } $, $ \mathbb{V}_{\left[ m,n,p\right] } $, $ \mathbb{I}_{\left[ m,n,p\right] } $ as $ \mathbb{V}_1 $,$ \mathbb{V}_2 $,$ \mathbb{V}_3 $,$ \mathbb{I} $, respectively. 
 For the sake of convenience, we will represent the variables $\mathcal{X},\mathcal{Y}$, and $\mathcal{C}$ in tensor form as vectors, as follows:\par
 $$
 x=V_{1} \left ( \mathcal{X}  \right ) , y=V_{2} \left ( \mathcal{Y}  \right ), c=V_{3} \left ( \mathcal{C}  \right )
 $$
 \par Then, we rewrite the objective function (\ref{eq:2}) as
 \begin{equation}\label{obj}
\tilde{F} \left ( x,y,c \right ) =G\left ( x,y,c \right ) +H\left ( c \right ) +\Phi \left ( c \right ) 
 \end{equation}
where $ x\in\mathbb{R}^{d_1},  y\in\mathbb{R}^{d_2}, c\in\mathbb{R}^{d_3}$
\begin{align*}
& \tilde{F}\left( x,y,c\right) : = G\left( x,y,c\right) + H\left ( c \right ) +\Phi \left ( c \right ),\\
& G\left( x,y,c\right) =\frac{1}{2}\Vert\mathbb{V}_1^{-1}\left( x\right)*\mathbb{V}_2^{-1}\left( y\right)-\mathbb{V}_3^{-1}\left( c\right)\Vert_{\ell_2}^2,\\
& H\left( c\right) : = \alpha_1\Vert D_1\left( c\right) \Vert_{\ell_1}+\alpha_2\Vert D_2\left( c\right) \Vert_{\ell_1},\\
& D_1\left( c\right) : = \mathcal{D}_1*_v\left[ \mathbb{V}_3^{-1}\left( c\right) \right] ,\\
& D_2\left( c\right) : = \left[ \mathbb{V}_3^{-1}\left( c\right) \right]*_v\mathcal{D}_2 ,\\
& \tilde{S}: =\left\lbrace w\in\mathbb{R}^{d_3}|w\left( i\right) =\left[ \mathbb{V}_3\left( \mathcal{G}\right) \right] \left( i\right) \mbox{ for }i \in \mathbb{I}\left( \Omega\right),w\left( i\right) \in\left[ 0,1\right]  \mbox{ for } i \notin \mathbb{I}\left( \Omega\right) \right\rbrace ,\\
&\Phi \left ( c \right )=\begin{cases} 0,c\in \tilde{S}\\ +\infty , c\notin  \tilde{S}\end{cases} 
\end{align*}
\par Because of the bilinear operation of "$ *_v $," it is evident that both $ D_1:\mathbb{R}^{m n p\times 1}\rightarrow\mathbb{R}^{m\times n\times p} $ and $ D_2:\mathbb{R}^{m n p\times 1}\rightarrow\mathbb{R}^{m\times n\times p} $ are linear operators. Moreover, it is easy to see that $ \tilde{S} $ is a non-empty closed set, which means that $ \Phi \left ( c \right )  $ is a proper lower semi-continuous (PLSC) function on $ \mathbb{R}^{d_3} $ and $ \tilde{F}\left( \cdot\right) $ is a PLSC function on $ \mathbb{R}^d $.
(\ref{eq:33}) is equivalent to
 \begin{equation}\label{obj2}
\begin{cases} x^{k+1}\in\mathop{\arg\min}_{x}\tilde{F}\left( x,y^k,c^k\right)+\frac{\rho_1}{2}\Vert x-x^k\Vert_2^2,\\ y^{k+1}\in\mathop{\arg\min}_{y}\tilde{F}\left( x^{k+1},y,c^k\right)+\frac{\rho_2}{2}\Vert y-y^k\Vert_2^2, \\ c^{k+1}\in\mathop{\arg\min}_{c}\tilde{F}\left( x^{k+1},y^{k+1},c\right)+\frac{\rho_3}{2}\Vert c-c^k\Vert_2^2 \end{cases}
 \end{equation}

 \par To show the convergence of PAM algorithm, the following convergence theory is needed.
 \begin{definition}{\rm{}(K{\L} property{\rm\cite{attouch2013convergence})}}
 	\begin{enumerate}[(a)]
 		\item The function $ f:\mathbb{R}^n\rightarrow\mathbb{R}\cup\left\lbrace +\infty\right\rbrace  $ is said to have the K{\L} property at $ \bar{x}\in dom\left( \partial f\right)  $, if there exist $ \eta\in\left( 0,+\infty\right]  $, a neighbourhood U of $ \bar{x} $ and a continuous	concave function $ \phi:\left[ 0,\eta\right) \rightarrow\left[ 0,+\infty\right)  $, such that:
 		\begin{enumerate}[(i)]
 			\item $ \phi\left( 0\right) =0 $
 			\item $ \phi $ is first-order continuous on $ \left( 0,\eta\right)  $ ,
 			\item $ \phi ' $ is positive on $ \left( 0,\eta\right) $ ,
 			\item for each $ x\in U\cap\left[ f\left( \bar{x}\right)<f<f\left( \bar{x}\right)+\eta  \right]  $, the K{\L} inequality holds:
 			\begin{equation*}
 			\phi '\left( f\left( x\right)-f\left( \bar{x}\right)  \right)dist\left( 0,\partial f\left( x\right) \right) \geq 1.
 			\end{equation*}
 		\end{enumerate}
 		\item PLSC functions which satisfies the K{\L} property at each point of $ dom\left( \partial f\right)  $ are called K{\L} functions, where the norm involved in $ dist\left( \cdot,\cdot\right)  $ is $ \Vert\cdot\Vert_2 $ and the convention $ dist\left( 0,\emptyset\right) : = +\infty $.
 	\end{enumerate}
 \end{definition}
 
\begin{lemma}{\cite{attouch2013convergence}}\label{lem1}
	Let $ f:\mathbb{R}^n\rightarrow\mathbb{R}\cup\left\lbrace +\infty\right\rbrace $ be a PLSC function. Let $ \left\lbrace x^{k}\right\rbrace_{k\in\mathbb{N}}\subset\mathbb{R}^n  $ be a sequence such that
	\begin{itemize}
		\item[H1] $ \left( \mbox{Sufficient decrease condition}\right) $ For each $ k\in\mathbb{N} $, there exits $ a\in\left( 0,+\infty\right) $ such that $ f\left( x^{k+1}\right)+a\Vert x^{k+1}-x^{k}\Vert_2^2\leq f\left( x^{k}\right)  $ hold
		\item[H2] $ \left( \mbox{Relative error condition}\right) $ For each $ k\in\mathbb{N} $, there exits $  w^{k+1}\in\partial f\left( x^{k+1}\right) $ and a constant $ b\in\left( 0,+\infty\right)  $ such that$ \Vert w^{k+1}\Vert_2\leq b\Vert x^{k+1}-x^{k}\Vert_2 $ hold
		\item[H3] $ \left( \mbox{Continuity condition}\right) $There exists a subsequence $ \left\lbrace x^{k_j}\right\rbrace_{j\in\mathbb{N}} $ and $ \bar{x}\in\mathbb{R}^n $ such that
		\begin{equation*}
		x^{k_j}\rightarrow\bar{x}\mbox{ and }f\left( x^{k_j}\right)\rightarrow f\left( \bar{x}\right),j\rightarrow\infty.
		\end{equation*}
	\end{itemize}
	If $f$ has the K{\L} property at $ \bar{x} $, then
	\begin{enumerate}[(i)]
		\item $ x^{k}\rightarrow\bar{x} $
		\item $ \bar{x} $ is a critical point of $f$, i.e., $ 0\in\partial f\left( \bar{x}\right)  $;
		\item the sequence $ \left\lbrace x^{k}\right\rbrace_k\in\mathbb{N} $ has a finite length, i.e.,
		\begin{equation*}
		\sum_{k=0}^{+\infty}\Vert x^{k+1}-x^{k}\Vert_2<+\infty.
		\end{equation*}
	\end{enumerate}
\end{lemma}
\par Second, we prove that the objective function  $\tilde{F} $ in (\ref{obj}) and the iterative sequence $ ( x^{k},y^{k},c^{k})_{k\in\mathbb{N}}  $ generated by iteration (\ref{obj2}) satisfy the assumptions in LEMMA \ref{lem1}. Thus we establish the convergence of the proposed algorithm.

For convenience, we denote $ \hat{S}=\mathbb{R}^{d_1}\times\mathbb{R}^{d_2}\times\tilde{S}\subset \mathbb{R}^d$.
Next, we need to prove that $ \tilde{F} $ satisfies the K{\L} property at each $ \left ( x,y,c \right ) \in\hat{S} $.
\begin{lemma}\label{lem2}
	$ \tilde{F} $ satisfies the K{\L} property at each $ \left ( x,y,c \right ) \in\hat{S} $.
\end{lemma}
\begin{proof}
For any $\left ( x,y,c \right )\in\hat{S}$, $ \tilde{F} $ can be expressed as follows:
\begin{equation*}
\tilde{F}: = G\left( x,y,c\right) +H\left( c\right)
\end{equation*}
Obviously, $ \hat{S} $ is a semi-algebraic set(\cite{lin2020tensor} Proposition 2).\par
For convenience, we denote $ u: =\left( u_1,u_2,u_3\right) , u_1:=x, u_2:=y, u_3:=c$ . $ D_1 $ and $ D_2 $ are linear mappings between finite dimensional spaces, each element of $ D_1\left( c\right) $ and $ D_2\left( c\right) $ can be regarded as a linear polynomial of $ u:=(x,y,c) $, so $ H\left( c\right) $ is only a finite sum of the composition of the absolute value function and linear polynomial of $u$. Therefore, $H\left( c\right) $ is semi-algebraic on $ \hat{S} $(\cite{lin2020tensor} Proposition 2,3). $ \mathbb{V}_1^{-1}\left( x\right) $, $ \mathbb{V}_2^{-1}\left( y\right) $, and $ \mathbb{V}_3^{-1}\left( c\right) $ are linear mappings between finite-dimensional spaces. Therefore, each element of $ \mathbb{V}_i^{-1}\left( u_i\right) \left( i=1,2,3\right) $ is a linear polynomial of $ u: =\left( u_1,u_2,u_3\right) $. In addition, according to the definition of $ \Vert\cdot\Vert_{\ell_2} $, we know that $ G\left( u\right) $ is a polynomial of $ u $. Thus, $G$ is semi-algebraic on $ \hat{S} $(\cite{lin2020tensor} Proposition 2). \par
In a word, $ \tilde{F} $ is semi-algebraic on $ \hat{S} $(\cite{lin2020tensor} Proposition 2). A semi-algebraic function is a K{\L} function (see\cite{bolte2007lojasiewicz,bolte2007clarke}).
Hence, $ \tilde{F} $ satisfies the K{\L} property at each $ \left ( x,y,c \right )\in\hat{S} $.
\end{proof}

\begin{theorem}
	Assume that the sequence $ u^{k}=( x^{k},y^{k},c^{k}) $ generated by iteration \rm{(\ref{obj2})} is bounded, then it converges to a critical point of $ \tilde{F} $.
\end{theorem}

\begin{proof}
As mentioned above, $ \tilde{F} $ is a PLSC function on $ \mathbb{R}^d $. From (\ref{obj2}), and we see that
\begin{align*}
\begin{cases} \tilde{F}( x^{k+1},y^{k},c^{k}) +\frac{\rho_1}{2}\Vert x^{k+1}-x^{k}\Vert_2^2\leq\tilde{F}(x^{k},y^{k},c^{k}),k\in\mathbb{N}\\ \tilde{F}( x^{k+1},y^{k+1},c^{k}) +\frac{\rho_2}{2}\Vert y^{k+1}-y^{k}\Vert_2^2\leq\tilde{F}( x^{k+1},y^{k},c^{k}),k\in\mathbb{N} \\ \tilde{F}( x^{k+1},y^{k+1},c^{k+1}) +\frac{\rho_3}{2}\Vert c^{k+1}-c^{k}\Vert_2^2\leq\tilde{F}( x^{k+1},y^{k+1},c^{k}),k\in\mathbb{N}. \end{cases}
\end{align*}
Combining the above three inequalities, we get
\begin{equation}\label{eq5}
\tilde{F}( u^{k+1}) +\frac{\rho_{\rm min}}{2}\Vert u^{k+1}-u^{k}\Vert_2^2\leq\tilde{F}( u^{k}) ,k\in\mathbb{N}.
\end{equation}
where $\rho_{\rm min}=\min\{\rho_1,\rho_2,\rho_3\}$. 
Therefore, H1(sufficient decrease condition) is satisfied, and $ a=\rho_{\rm min}/2 $.\par
Denote $ u: =\left( u_1,u_2,u_3\right) , u_1:=x, u_2:=y, u_3:=c$, $ M=\lbrace u^{k}|k\in\mathbb{N} \rbrace $, for $ u=( u_1,u_2,u_3)\in\mathbb{R}^s $, $ u_i\in\mathbb{R}^{d_i}( i=1,2,3)  $.
 Because $G$ is a polynomial function, so it is infinitely differentiable. From (\ref{obj2}) and Proposition 1\cite{lin2020tensor},

\begin{equation}\label{obj5}
\begin{cases} 0\in \nabla_{u_1}G( u_1^{k+1},u_2^{k},u_3^{k})+\rho_1( u_1^{k+1}-u_1^{k}) ,k\in\mathbb{N}\\ 0\in\nabla_{u_2}G( u_1^{k+1},u_2^{k+1},u_3^{k})+\rho_2( u_2^{k+1}-u_2^{k}) ,k\in\mathbb{N}\\ 0\in\nabla_{u_3}G( u_1^{k+1},u_2^{k+1},u_3^{k+1})+\partial( \Phi+H)( u_3^{k+1}) +\rho_3( u_3^{k+1}-u_3^{k}) ,k\in\mathbb{N}. \end{cases}
\end{equation}

(\ref{obj5}) means the existence of $ w_1^{k+1}\in\partial( \Phi+H)( u_3^{k+1}) $, so that
\begin{equation*}			-\tilde{w}^{k+1}=\tilde{u}^{k}+(\rho_1(u_1^{k+1}-u_1^{k}),\rho_2(u_2^{k+1}-u_2^{k}),\rho_3(u_3^{k+1}-u_3^{k})) ,k\in\mathbb{N},
\end{equation*}
where
\begin{align*}
&\tilde{w}^{k+1}=( 0,0,w_1^{k+1}) ,\\
&\tilde{u}^{k}=( \nabla_{u_1}G( u_1^{k+1},u_2^{k},u_3^{k}),\nabla_{u_2}G( u_1^{k+1},u_2^{k+1},u_3^{k}),\nabla_{u_3}G( u_1^{k+1},u_2^{k+1},u_3^{k+1}) ).
\end{align*}
\par Let $ w^{k+1}=\tilde{w}^{k+1}+\nabla G( u^{k+1}) $, obviously $ w^{k+1}\in\partial\tilde{F}( u^{k+1}) $ and
\begin{equation}\label{eq9}
-w^{k+1}= \tilde{u}^{k}-\nabla G( u^{k+1}) +(\rho_1(u_1^{k+1}-u_1^{k}),\rho_2(u_2^{k+1}-u_2^{k}),\rho_3(u_3^{k+1}-u_3^{k})) ,k\in\mathbb{N},
\end{equation}
\par  Define the coordinate projections by
\begin{equation*}
\Psi_i( u) : = u_i,i=1,2,3.
\end{equation*}
\par Denote $ \hat{M}=\Psi_1( u)\times\Psi_2( u)\times\Psi_3( u) $, since $M$ is bounded, so is $ \hat{M}\subset \mathbb{R}^d $. Since $G$ is a polynomial, it is easy to prove that $ \nabla G $ is Lipschitz continuous on any bounded subset of $ \mathbb{R}^d  $. Therefore, for any $ u,w\in\hat{M} $, there is a constant $ c>0 $ such that
\begin{equation*}
\Vert\nabla G( u) -\nabla G( w) \Vert_2^2\leq c\Vert u-w\Vert_2^2
\end{equation*}
hence,
\begin{align*}
\Vert \tilde{u}^{k}-\nabla G( u^{k+1}) \Vert_2=&\left[ \Vert\nabla_{u_1}G( u_1^{k+1},u_2^{k},u_3^{k})-\nabla_{u_1}G( u^{k+1}) \Vert_2^2 \right.\\ &\left. +\Vert\nabla_{u_2}G( u_1^{k+1},u_2^{k+1},u_3^{k})-\nabla_{u_2}G( u^{k+1}) \Vert_2^2\right] ^{\frac{1}{2}}\\
\leq & \left[ c^2( \Vert u_2^{k}-u_2^{k+1}\Vert_2^2+\Vert u_3^{k}-u_3^{k+1}\Vert_2^2)+c^2( \Vert u_3^{k}-u_3^{k+1}\Vert_2^2)  \right] ^{\frac{1}{2}}\\
\leq & \sqrt{2}c\Vert u^{k+1}-u^{k}\Vert_2^2,k\in\mathbb{N}
\end{align*}
according to (\ref{eq9}), there are
\begin{align*}
\Vert w^{k+1}\Vert_2 & \leq\Vert \tilde{u}^{k}-\nabla G( u^{(k+1)}) \Vert_2+\rho_{\max}\Vert u^{k+1}-u^{k}\Vert_2\\
& \leq ( \sqrt{2}c+\rho_{\max}) \Vert u^{k+1}-u^{k}\Vert_2,k\in\mathbb{N}
\end{align*}
where $\rho_{\rm{max}}=\max\{\rho_1,\rho_2,\rho_3\}$. Therefore, H2(Relative error condition) is satisfied, and $ b=\sqrt{2}c+\rho_{\max} $.\par
In addition, because $ \left\lbrace u^{k}\right\rbrace_{k\in\mathbb{N}}\subset\mathbb{R}^d  $ is bounded, it is relatively compact, so there exists a subsequence $ \left\lbrace u^{k_j}\right\rbrace _{j\in\mathbb{N}} $ and $ \bar{u}\in\mathbb{R}^d $ such that $ u^{k_j}\rightarrow\bar{u} $, $ j\rightarrow +\infty $. Since $ \left\lbrace u_3^{k}|k\in\mathbb{N}\right\rbrace\subset\tilde{S}  $, so $ \left\lbrace u^{k}|k\in\mathbb{N}\right\rbrace\subset\hat{S}  $ holds. Since $ \hat{S} $ is closed, $ \bar{u}\in\hat{S} $, and $ \tilde{F} $ is continuous on $ \hat{S} $. Therefore, $ \tilde{F}\left( u^{k_j}\right)\rightarrow\tilde{F}\left( \bar{u}\right)   $ , $ j\rightarrow +\infty $. Hence, H3(Continuity condition) is satisfied.\par
Finally, according to Lemma \ref{lem2}, $ \tilde{F} $ satisfies the K{\L} property at $ \bar{u}\in\hat{S} $. By Lemma \ref{lem1}, the sequence $\{u^{k}\}_{k\in \mathbb{N}}$ converges to the critical point of $ \tilde{F} $. This completes the proof.
\end{proof}


\bibliographystyle{siamplain}

\begin{thebibliography}{10}

\bibitem{attouch2013convergence}
{\sc H.~Attouch, J.~Bolte, and B.~F. Svaiter}, {\em Convergence of descent
  methods for semi-algebraic and tame problems: proximal algorithms,
  forward--backward splitting, and regularized gauss--seidel methods},
  Mathematical Programming, 137 (2013), pp.~91--129.

\bibitem{bolte2007lojasiewicz}
{\sc J.~Bolte, A.~Daniilidis, and A.~Lewis}, {\em The {\l}ojasiewicz inequality
  for nonsmooth subanalytic functions with applications to subgradient
  dynamical systems}, SIAM Journal on Optimization, 17 (2007), pp.~1205--1223.

\bibitem{bolte2007clarke}
{\sc J.~Bolte, A.~Daniilidis, A.~Lewis, and M.~Shiota}, {\em Clarke
  subgradients of stratifiable functions}, SIAM Journal on Optimization, 18
  (2007), pp.~556--572.

\bibitem{Br10}
{\sc K.~Braman}, {\em Third-order tensors as linear operators on a space of
  matrices}, Linear Algebra and Its Applications, 433 (2010), pp.~1241--1253.

\bibitem{chen2015fractional}
{\sc D.~Chen, Y.~Chen, and D.~Xue}, {\em Fractional-order total variation image
  denoising based on proximity algorithm}, Applied Mathematics and Computation,
  257 (2015), pp.~537--545.

\bibitem{chen2022color}
{\sc J.~Chen and M.~K. Ng}, {\em Color image inpainting via robust pure
  quaternion matrix completion: Error bound and weighted loss}, SIAM Journal on
  Imaging Sciences, 15 (2022), pp.~1469--1498.

\bibitem{CXZ20}
{\sc Y.~Chen, X.~Xiao, and Y.~Zhou}, {\em Multi-view subspace clustering via
  simultaneously learning the representation tensor and affinity matrix},
  Pattern Recognition, 106 (2020).

\bibitem{deng2023t}
{\sc Y.-J. Deng, H.-C. Li, S.-Q. Tan, J.~Hou, Q.~Du, and A.~Plaza}, {\em
  t-linear tensor subspace learning for robust feature extraction of
  hyperspectral images}, IEEE Transactions on Geoscience and Remote Sensing, 61
  (2023), pp.~1--15.

\bibitem{DT15}
{\sc C.~Donciu and M.~Temneanu}, {\em An alternative method to zero-padded
  {DFT}}, Measurement, 70 (2015), pp.~14--20.

\bibitem{GV13}
{\sc G.~H. Golub and C.~F. Van~Loan}, {\em Matrix Computation 4th ed.}, The
  Johns Hopkins University Press, Baltimore, USA, 2013.

\bibitem{HZW21}
{\sc J.~Hou, F.~Zhang, and J.~Wang}, {\em One-bit tensor completion via
  transformed tensor singular value decomposition}, Applied Mathematical
  Modelling, 95 (2021), pp.~760--782.

\bibitem{JNZH20}
{\sc T.~X. Jiang, M.~K. Ng, X.~L. Zhao, and T.~Z. Huang}, {\em Framelet
  representation of tensor nuclear norm for third-order tensor completion},
  IEEE Transaction on Image Processing, 29 (2020), pp.~7233--7244.

\bibitem{KBHH13}
{\sc M.~E. Kilmer, K.~Braman, N.~Hao, and R.~Hoover}, {\em Third-order tensors
  as operators on matrices: A theoretical and computational framework with
  applications in imaging}, SIAM Journal on Matrix Analysis and Applications,
  34 (2013), pp.~148--172.

\bibitem{KM11}
{\sc M.~E. Kilmer and C.~D. Martin}, {\em Factorization strategies for
  third-order tensors}, Linear Algebra and Its Applications, 435 (2011),
  pp.~641--658.

\bibitem{KKA15}
{\sc E.~Knerfeld, M.~E. Kilmer, and S.~Aeron}, {\em Tensor-tensor products with
  invertible linear transforms}, Linear Algebra and Its Applications, 485
  (2015), pp.~545--570.

\bibitem{kolda2009tensor}
{\sc T.~G. Kolda and B.~W. Bader}, {\em Tensor decompositions and
  applications}, SIAM review, 51 (2009), pp.~455--500.

\bibitem{LS03}
{\sc N.~Le~Bihan and S.~Sangwine}, {\em Quaternion principal component analysis
  of color images}, IEEE International Conference on Image Processing, 1
  (2003), pp.~809--812.

\bibitem{li2022nonlinear}
{\sc B.-Z. Li, X.-L. Zhao, T.-Y. Ji, X.-J. Zhang, and T.-Z. Huang}, {\em
  Nonlinear transform induced tensor nuclear norm for tensor completion},
  Journal of Scientific Computing, 92 (2022), p.~83.

\bibitem{lin2020tensor}
{\sc X.-L. Lin, M.~K. Ng, and X.-L. Zhao}, {\em Tensor factorization with total
  variation and tikhonov regularization for low-rank tensor completion in
  imaging data}, Journal of Mathematical Imaging and Vision, 62 (2020),
  pp.~900--918.

\bibitem{LLOQ21}
{\sc C.~Ling, J.~Liu, C.~Ouyang, and L.~Qi}, {\em St-svd factorization and
  s-diagonal tensors}, arXiv:2104.05329,  (2021).

\bibitem{liu2014generalized}
{\sc R.~W. Liu, L.~Shi, W.~Huang, J.~Xu, S.~C.~H. Yu, and D.~Wang}, {\em
  Generalized total variation-based mri rician denoising model with spatially
  adaptive regularization parameters}, Magnetic Resonance Imaging, 32 (2014),
  pp.~702--720.

\bibitem{MKL20}
{\sc J.~Miao, K.~I. Kou, and W.~Liu}, {\em Low-rank quaternion tensor
  completion for recovering color videos and images}, Pattern Recognition, 107
  (2020).

\bibitem{MWGdD02}
{\sc B.~Muquet, Z.~Wang, G.~B. Giannakis, M.~de~Courville, and P.~Duhamel},
  {\em Cyclic prefixing or zero padding for wireless multicarrier
  transmissions?}, IEEE Transactions on Communications, 5 (2002),
  pp.~2136--2148.

\bibitem{QLLO21}
{\sc L.~Qi, C.~Ling, J.~Liu, and C.~Ouyang}, {\em An orthogonal equivalence
  theorem for third order tensors}, to appear in: Journal of Industrial and
  Management Optimization,  (2021).

\bibitem{QL21}
{\sc L.~Qi and Z.~Luo}, {\em Tubal matrices}, arXiv:2105.00793v3,  (2021).

\bibitem{qi2022quaternion}
{\sc L.~Qi, Z.~Luo, Q.-W. Wang, and X.~Zhang}, {\em Quaternion matrix
  optimization: Motivation and analysis}, Journal of Optimization Theory and
  Applications, 193 (2022), pp.~621--648.

\bibitem{qiu2021nonlocal}
{\sc D.~Qiu, M.~Bai, M.~K. Ng, and X.~Zhang}, {\em Nonlocal robust tensor
  recovery with nonconvex regularization}, Inverse Problems, 37 (2021),
  p.~035001.

\bibitem{qiu2021robust}
{\sc D.~Qiu, M.~Bai, M.~K. Ng, and X.~Zhang}, {\em Robust low transformed
  multi-rank tensor methods for image alignment}, Journal of Scientific
  Computing, 87 (2021), p.~24.

\bibitem{SHKM14}
{\sc O.~Semerci, N.~Hao, M.~E. Kilmer, and E.~L. Miller}, {\em Tensor-based
  formulation and nuclear norm regularization for multienergy computed
  tomography}, IEEE Transactions on Image Processing, 23 (2014),
  pp.~1678--1693.

\bibitem{shi2021robust}
{\sc Q.~Shi, Y.-M. Cheung, and J.~Lou}, {\em Robust tensor svd and recovery
  with rank estimation}, IEEE Transactions on Cybernetics, 52 (2021),
  pp.~10667--10682.

\bibitem{SNZ20}
{\sc G.~Song, M.~K. Ng, and X.~Zhang}, {\em Robust tensor completion using
  transformed tensor singular value decomposition}, Numerical Linear Algebra
  with Applications, 27 (2020).

\bibitem{tarzanagh2018fast}
{\sc D.~A. Tarzanagh and G.~Michailidis}, {\em Fast randomized algorithms for
  t-product based tensor operations and decompositions with applications to
  imaging data}, SIAM Journal on Imaging Sciences, 11 (2018), pp.~2629--2664.

\bibitem{WGLZ21}
{\sc X.~Wang, L.~Gu, H.~W. Lee, and G.~Zhang}, {\em Quantumn context-aware
  recommendation systems based on tensor singular value decomposition}, Quantum
  Information Processing, 20 (2021).

\bibitem{xiao2020prior}
{\sc X.~Xiao, Y.~Chen, Y.-J. Gong, and Y.~Zhou}, {\em Prior knowledge
  regularized multiview self-representation and its applications}, IEEE
  Transactions on neural networks and learning systems, 32 (2020),
  pp.~1325--1338.

\bibitem{XCGZ21}
{\sc X.~Xiao, Y.~Chen, Y.~J. Gong, and Y.~Zhou}, {\em Low-rank reserving
  t-linear projection for robust image feature extraction}, IEEE Transactions
  on Image Processing, 30 (2021), pp.~108--120.

\bibitem{xu2013parallel}
{\sc Y.~Xu, R.~Hao, W.~Yin, and Z.~Su}, {\em Parallel matrix factorization for
  low-rank tensor completion}, arXiv preprint arXiv:1312.1254,  (2013).

\bibitem{YHHH16}
{\sc L.~Yang, Z.~H. Huang, S.~Hu, and J.~Han}, {\em An iterative algorithm for
  third-order tensor multi-rank minimization}, Computational Optimization and
  Applications, 63 (2016), pp.~169--202.

\bibitem{yu2022t}
{\sc Q.~Yu and X.~Zhang}, {\em T-product factorization based method for matrix
  and tensor completion problems}, Computational Optimization and Applications,
   (2022), pp.~1--28.

\bibitem{yu2020multi}
{\sc Q.~Yu, X.~Zhang, and Z.-H. Huang}, {\em Multi-tubal rank of third order
  tensor and related low rank tensor completion problem}, arXiv preprint
  arXiv:2012.05065,  (2020).

\bibitem{ZSKA18}
{\sc J.~Zhang, A.~K. Saibaba, M.~E. Kilmer, and S.~Aeron}, {\em A randomized
  tensor singular value decomposition based on the t-product}, Numerical Linear
  Algebra with Applications, 25 (2018).

\bibitem{ZA17}
{\sc Z.~Zhang and S.~Aeron}, {\em Exact tensor completion using t-{SVD}}, IEEE
  Transactions on Signal Processing, 65 (2017), pp.~1511--1526.

\bibitem{ZEAHK14}
{\sc Z.~Zhang, G.~Ely, S.~Aeron, N.~Hao, and M.~E. Kilmer}, {\em Novel methods
  for multilinear data completion and de-noising based on tensor-{SVD}}, in
  Proceedings of the IEEE Conference on Computer Vision and Pattern
  Recognition, ser. CVPR '14, IEEE, 2014, pp.~3842--3849.

\bibitem{zhao2022robust}
{\sc X.~Zhao, M.~Bai, D.~Sun, and L.~Zheng}, {\em Robust tensor completion:
  Equivalent surrogates, error bounds, and algorithms}, SIAM Journal on Imaging
  Sciences, 15 (2022), pp.~625--669.

\bibitem{zheng2020tensor}
{\sc Y.-B. Zheng, T.-Z. Huang, X.-L. Zhao, T.-X. Jiang, T.-Y. Ji, and T.-H.
  Ma}, {\em Tensor n-tubal rank and its convex relaxation for low-rank tensor
  recovery}, Information Sciences, 532 (2020), pp.~170--189.

\bibitem{zheng2019mixed}
{\sc Y.-B. Zheng, T.-Z. Huang, X.-L. Zhao, T.-X. Jiang, T.-H. Ma, and T.-Y.
  Ji}, {\em Mixed noise removal in hyperspectral image via low-fibered-rank
  regularization}, IEEE Transactions on Geoscience and Remote Sensing, 58
  (2019), pp.~734--749.

\bibitem{ZLLZ18}
{\sc P.~Zhou, C.~Lu, Z.~Lin, and C.~Zhang}, {\em Tensor factorization for
  low-rank tensor completion}, IEEE Transactions on Image Processing, 27
  (2018), pp.~1152--1163.

\bibitem{ZC21}
{\sc Y.~Zhou and Y.~Cheung}, {\em Bayesian low-tubal-rank robust tensor
  factorization with multi-rank determination}, IEEE Transactions on Pattern
  Analysis and Machine Intelligence, 43 (2021), pp.~62--76.

\end{thebibliography}

\end{document}